\documentclass{article}
 \usepackage{amsthm,amsmath,amssymb,amsfonts, url}
\usepackage{graphicx}
\usepackage[colorlinks,urlcolor=blue]{hyperref}
\usepackage{epsfig}
\usepackage{kotex}
\theoremstyle{definition}
\newtheorem{definition}{Definition}[section]
\newtheorem{example}[definition]{Example}

\theoremstyle{plain}
\newtheorem{theorem}[definition]{Theorem}

\newtheorem{proposition}[definition]{Proposition}
\newtheorem{corollary}[definition]{Corollary}

\numberwithin{equation}{section}
\begin{document}

%%%%%%%%%%%%%%%%%%%%%%%%%%%%%%%%%%%%%%%%%%%%%%%%%%%%%%%%%%%%%%%%%%%%%%
\noindent {\large \textbf
	{Ordered homomorphisms and kernels of ordered BCI-algebras}}\\[13pt]
%%%%%%%%%%%%%%%%%%%%%%%%%%%%%%%%%%%%%%%%%%%%%%%%%%%%%%%%%%%%%%%%%%%%%%
{Eunsuk Yang$^{1, *}$,  Eun Hwan Roh$^{2}$ and Young Bae Jun$^{3}$} \\[1mm]
%%%%%%%%%%%%%%%%%%%%%%%%%%%%%%%%%%%%%%%%%%%%%%%%%%%%%%%%%%%%%%%%%%%%%%
$^{1}$Department of Philosophy, Jeonbuk National University, Jeonju 54896, Korea \\
e-mail: eunsyang@jbnu.ac.kr \\ [2mm]
%%%%%%%%%%%%%%%%%%%%%%%%%%%%%%%%%%%%%%%%%%%%%%%%%%%%%%%%%%%%%%%%%%%%%%
$^{2}$Department of Mathematics Education,
Chinju National University of Education, Jinju 52673, Korea \\
e-mail: ehroh9988@gmail.com \\ [2mm]
%%%%%%%%%%%%%%%%%%%%%%%%%%%%%%%%%%%%%%%%%%%%%%%%%%%%%%%%%%%%%%%%%%%%%%
$^{3}$Department of Mathematics Education,
Gyeongsang National University,  Jinju 52828, Korea \\
e-mail: skywine@gmail.com  \\ [2mm]
%%%%%%%%%%%%%%%%%%%%%%%%%%%%%%%%%%%%%%%%%%%%%%%%%%%%%%%%%%%%%%%%%%%%%%
*Corresponding Author: E. Yang (eunsyang@jbnu.ac.kr)

\noindent\hrulefill \\
%%%%%%%%%%%%%%%%%%%%%%%%%%%%%%%%%%%%%%%%%%%%%%%%%%%%%%%%%%%%%%%%%%%%%%
{\bf Abstract}
Recently Yang-Roh-Jun introduced the notion of ordered BCI-algebras as a generalization of BCI-algebras. They also introduced the notions of homomorphisms and kernels of ordered BCI-algebras and investigated related properties. Here we extend their investigation to ordered homomorphisms, i.e., order-preserving homomorphisms. To this end, the notions of ordered homomorphism and kernel of ordered BCI-algebras are first defined. Next, properties associated with (ordered) subalgebras, (ordered) filters and direct products of ordered BCI-algebras are addressed.
%{\tiny{(007-221215R07-Ordered homomorphisms and kernels of ordered BCI-algebras.tex)}}\\[1mm]
%%%%%%%%%%%%%%%%%%%%%%%%%%%%%%%%%%%%%%%%%%%%%%%%%%%%%%%%%%%%%%%%%%%%%%
%한글을 사용하려면 \verb|\usepackage{kotex}|을 사용해야 합니다.\\
%%%%%%%%%%%%%%%%%%%%%%%%%%%%%%%%%%%%%%%%%%%%%%%%%%%%%%%%%%%%%%%%%%%%%%

\vspace{2mm} \noindent 

\vspace{2mm} \noindent
%% 06F35 BCK-algebras, BCI-algebras (aspects of ordered structures) [See also 03G25]
%% 03G25 Other algebras related to logic [See also 03F45, 06D20, 06E25, 06F35]
%% 06D20 Heyting algebras (lattice-theoretic aspects) [See also 03G25]
%% 08A72 Fuzzy algebraic structures
%% 03B05 Classical propositional logic (Sheffer stroke)

\noindent\hrulefill \\
%%%%%%%%%%%%%%%%%%%%%%%%%%%%%%%%%%%%%%%%%%%%%%%%%%%%%%%%%%%%%%%%%%%%%%
%% \renewcommand{\thefootnote}{}
%% \footnotetext{*Corresponding author.}
%% \footnotetext{e-mail:  }
%%%%%%%%%%%%%%%%%%%%%%%%%%%%%%%%%%%%%%%%%%%%%%%%%%%%%%%%%%%%%%%%%%%%%%
\def\LomX{X}  \def\LomY{Y}  
\def\LomF{F}   \def\LomG{G}
\def\LomA{A}   \def\LomB{B}
 \def\op{\rightarrow}  \def\oq{\Rightarrow} 
\def\jel{\le_{\jI}}
%% \def\orl{\preceq} \def\org{\succeq}
%%  \def\Xo0{\textbf{\LomX} := (\LomX, \op, \jI, \le)}
%%%%%%%%%%%%%%%%%%%%%%%%%%%%%%%%%%%%%%%%%%%%%%%%%%%%%%%%%%%%%%%%%%%%%
\def\opX{\op_{\LomX}}  \def\opY{\op_{\LomY}} 
\def\leX{\le_{\LomX}}  \def\leY{\le_{\LomY}} 
\def\jIX{\jI_{\LomX}}  \def\jIY{\jI_{\LomY}} 
 \def\Xo0{$\textbf{\LomX} := (\LomX,$ $\op,$ $\jI,$ $\jel)$}
\def\Yo0{$\textbf{\LomY} := (\LomY,$ $\op,$ $\jI,$ $\jel)$}
%%%%%%%%%%%%%%%%%%%%%%%%%%%%%%%%%%%%%%%%%%%%%%%%%%%%%%%%%%%%%%%%%%%%%
%% \def\Ao0{$\textbf{\LomA} := (\LomA,$ $\op,$ $\jI,$ $\jel)$}
%%%%%%%%%%%%%%%%%%%%%%%%%%%%%%%%%%%%%%%%%%%%%%%%%%%%%%%%%%%%%%%%%%%%%%%
%%   \def\ja{a} \def\jb{b}  \def\jc{c}  \def\jd{d} %% \def\je{e} \def\jf{f}
%%%%%%%%%%%%%%%%%%%%%%%%%%%%%%%%%%%%%%%%%%%%%%%%%%%%%%%%%%%%%%%%%%%%%%
\def\0{0} \def\1{1} \def\2{\tfrac{3}{4}} \def\3{\tfrac{1}{2}} \def\4{\tfrac{1}{4}}  
\def\8{\tfrac{2}{3}} \def\9{\tfrac{1}{3}} 
\def\x{x}  \def\y{y} \def\z{z}
\def\jI{e} \def\jx{a}  \def\jy{b}  \def\jz{c}  \def\jd{d}
%%%%%%%%%%%%%%%%%%%%%%%%%%%%%%%%%%%%%%%%%%%%%%%%%%%%%%%%%%%%%%%%%%%%%%
\def\kap{f}   \def\kbp{g} \def\kcp{h}
%% \def\kap{\mu}   \def\kbp{\gamma}  \def\jkap{f} \def\jkbp{g}
%% \def\kapA{\kap_{A}}   \def\kbpA{\kbp_{A}}
%% \def\kapB{\kap_{B}}   \def\kbpB{\kbp_{B}}
%% \def\jn{\varepsilon}   \def\jm{\delta}
%% \def\jt{t}   \def\js{s}
%% \def\jle{\jle}
%% \def\Lkap{\L{}_{\kap}}   \def\Lkbp{\L{}_{\kbp}}
%%%%%%%%%%%%%%%%%%%%%%%%%%%%%%%%%%%%%%%%%%%%%%%%%%%%%%%%%%%%%%%%%%%%%%
%% \def\ivq{\in\! \vee  q\,}   %%  \def\iwq{\in\! \wedge  q\,}
%% \def\PI{positive implicative }
%% \def\la{[}\def\ra{]}
%% \def\la{\langle}\def\ra{\rangle}
%% \def\IF{intuitionistic fuzzy }
%%%%%%%%%%%%%%%%%%%%%%%%%%%%%%%%%%%%%%%%%%%%%%%%%%%%%%%%%%%%%%%%%%%%%%
%%%%%%%%%%%%%%%%%%%%%%%%%%%%%%%%%%%%%%%%%%%%%%%%%%%%%%%%%%%%%%%%%%%%%%
\def\jop{\circ}
\def\RXo0{$\textbf{\LomX} := (\LomX,$ $\op,$ $\jop,$ $\jI,$ $\jel)$}
\def\RYo0{$\textbf{\LomY} := (\LomY,$ $\op,$ $\jop,$ $\jI,$ $\jel)$}
%%%%%%%%%%%%%%%%%%%%%%%%%%%%%%%%%%%%%%%%%%%%%%%%%%%%%%%%%%%%%%%%%%%%%%
%%%%%%%%%%%%%%%%%%%%%%%%%%%%%%%%%%%%%%%%%%%%%%%%%%%%%%%%%%%%%%%%%%%%%%
\def\p{\partial}
\def\PXo0{$\textbf{\LomX} := (\LomX,$ $\op,$ $\jI,$ $\p,$ $\jel)$}
\def\joq{\sim}
\def\RPXo0{$\textbf{\LomX} := (\LomX,$ $\op,$ $\jop,$ $\jI,$ $\p,$ $\jel)$}
\def\RPYo0{$\textbf{\LomX} := (\LomY,$ $\op,$ $\jop,$ $\jI,$ $\p,$ $\jel)$}
%\def\ejl{='}
%\def\ree{=_{e}}
%%%%%%%%%%%%%%%%%%%%%%%%%%%%%%%%%%%%%%%%%%%%%%%%%%%%%%%%%%%%%%%%%%%%%%
%\def\1kap{\kap_{1}}  \def\opY{\op_{\LomY}} 
%\def\leX{\le_{\LomX}}  \def\leY{\le_{\LomY}} 
%\def\jIX{\jI_{\LomX}}  \def\jIY{\jI_{\LomY}} 
% \def\Xo0{$\textbf{\LomX} := (\LomX,$ $\op,$ $\jI,$ $\jel)$}
%\def\Yo0{$\textbf{\LomY} := (\LomY,$ $\op,$ $\jI,$ $\jel)$}
%%%%%%%%%%%%%%%%%%%%%%%%%%%%%%%%%%%%%%%%%%%%%%%%%%%%%%%%%%%%%%%%%%%%%
%\def\RXo0{$\textbf{\LomX} := (\LomX,$ $\op,$ $\jop,$ $\jI,$ $\jel)$}
%\def\RYo0{$\textbf{\LomY} := (\LomY,$ $\op,$ $\jop,$ $\jI,$ $\jel)$}
%%%%%%%%%%%%%%%%%%%%%%%%%%%%%%%%%%%%%%%%%%%%%%%%%%%%%%%%%%%%%%%%%%%%%%
%\def\loA{\LomX_{1}}  \def\loB{\LomX_{2}} 
%\def\loC{\LomY_{1}}  \def\loD{\LomY_{2}} 
%\def\lomZ{\LomX_{1}}  \def\lomW{\LomX_{2}} 
%\def\lomK{\LomY_{1}}  \def\lomL{\LomY_{2}} 
\def\LomK{K}
\def\jIA{\jI_{\LomX1}}  \def\jIB{\jI_{\LomX2}} 
\def\jIC{\jI_{\LomY1}}  \def\jID{\jI_{\LomY2}} 
\def\opA{\op_{\LomX1}}  \def\opB{\op_{\LomX2}} 
\def\opC{\op_{\LomY1}}  \def\opD{\op_{\LomY2}} 
\def\leA{\le_{\LomX1}}  \def\leB{\le_{\LomX2}} 
\def\leC{\le_{\LomY1}}  \def\leD{\le_{\LomY2}} 
\def\oqX{\oq_{\LomX}}  \def\llX{\ll_{\LomX}} 
\def\oqY{\oq_{\LomY}}  \def\llY{\ll_{\LomY}} 
\def\ix{x_{1}} \def\kx{x_{2}} \def\iy{y_{1}} \def\ky{y_{2}}
\def\jxi{\jx_{1}} \def\jxk{\jx_{2}} \def\jyi{\jy_{1}} \def\jyk{\jy_{2}}
%\def\aX{$\textbf{\LomX}_{1}$ := ($\LomX_{1}$, $\opX$$_{1}$, $\jIX$$_{1}$, $\leX$$_{1})$
%\def\XDo0{$\textbf{\LomX}_{1}$ := $(\LomX_{1}, \opA$, $\jIA$, $\leA)$}
%%%%%%%%%%%%%%%%%%%%%%%%%%%%%%%%%%%%%%%%%%%%%%%%%%%%%%%%%%%%%%%%%%%%%%
%% \linenumbers  \def\contentsname{Contents}
%\tableofcontents
%% \thispagestyle{empty}   \linenumbers

%%%%%%%%%%%%%%%%%%%%%%%%%%%%%%%%%%%%%%%%%%%%%%%%%%%%%%%%%%%%%%%%%%%%%%
\section{Introduction: Background and Motivation} 
%\section{Introduction and Preliminaries} 
%%%%%%%%%%%%%%%%%%%%%%%%%%%%%%%%%%%%%%%%%%%%%%%%%%%%%%%%%%%%%%%%%%%%%%

One of important rsearch areas of universal algebra and logic is to introduce algebras and logics with more general structures and to study related properties. After Imai and Is\'eki \cite{ImIs66} introduced the notion of BCK-algebras as logic algebras,  a lot of algebras and logics have been introduced following this research trend. For instance, in the same year Is\'eki \cite{Is66} introduced the notion of BCI-algebras as a generalization of BCK-algebras. Since then, many algebras such as BCH-algebras \cite{HL83, HL85}, BH-algebras \cite{JRK98} and B-algebras \cite{NK02a, NK02b} have introduced as generalizations of BCI-algebras. Moreover, homomorphisms of such algebras have been studied (see e.g. \cite{H-Book, AT89, DGW08, DJ02, CF03, SC16, JRK98, JKK03, NN18, NK02b, AL10, SV20}).  

Very recently, Yang-Roh-Jun \cite{221005} introduced the notion of ordered BCI-algebras (briefly OBCI-algebras) as a generalization of BCI-algebras. They \cite{221215} further introduced the notions of homomorphisms of OBCI-algebras and studied associated properties. Here we note that ordered algebras were introduced as an algebraic generalization with a partial order. More precisely, ordered algebras were defined as structures with an underlying partial order and operations preserving orders, i.e., monotone operations (see e.g. \cite{B76,DH01,F63}). Although OBCI-algebras are a generalization of BCI-algebras with an underlying partial order, they are not ordered algebras in the above sense. As \eqref{(a3)} below shows, an implication operation in an OBCI-algebra is not mononotone with respect to its consequent argument. Because of such features of operations, Dunn \cite{D93} introduced tonoids as ordered sets with operations, each argument of which is isotone or antitone.  

Note that ordered algebras have been recently defined in more general contexts. For instance, Jansana-Moraschini \cite{JM18} defined the notion of ordered algebras in the sense of the notion of Dunn's tonoids. Raftery \cite{APAL164-251} defined partially ordered algebras as structures consisting of algebras and partial orders. In a similar context, Yang-Dunn \cite{YD21a} introduced implicational tonoids as implicational ordered structures. All these structures do not require the condition that operations must be monotone. In this sense, the OBCI-algebras can be regarded as ordered algebras considered in a more general context, i.e., ordered algebras dropping the monotonicity condition of operations.   

One interesting fact is that homomorphisms of BCI-algebras are order-preserving in themselves (see \cite{DGW08,H-Book}), whereas homomorphisms of BCI-algebras need not (see Example \ref{ETkcGV70-230520} below). This means that the order-preservation is not a property of homomorphisms of OBCI-algebras. This situation provides a motivation for us to study a mapping of OBCI-algebras that satisfies both homomorphism and monotinicity. To this end, we introduce the notion of ordered homomorphisms of OBCI-algebras as order-preserving homomorphisms of OBCI-algebras and investigate related properties. More exactly, we first define ordered homomorphisms and kernels of OBCI-algebras. We then address properties related to (ordered) subalgebras, (ordered) filters and direct products of OBCI-algebras in order to emphasize the similarities and differences with the research of homomorphisms of OBCI-algebras in \cite{221215}.

%%%%%%%%%%%%%%%%%%%%%%%%%%%%%%%%%%%%%%%%%%%%%%%%%%%%%%%%%%%%%%%%%%%%%%
\section{Preliminaries} 
%%%%%%%%%%%%%%%%%%%%%%%%%%%%%%%%%%%%%%%%%%%%%%%%%%%%%%%%%%%%%%%%%%%%%%

\begin{definition}[\cite{221005}] \label{def31-221006}	
Let $\LomX$ be a set with a binary operation $``\op"$,  a constant $``\jI"$ and a binary relation $``\jel"$. Then 
\Xo0 is called an {\it ordered BCI-algebra}
(briefly, OBCI-algebra) if it satisfies the following conditions: 
\begin{align}
 &\label{OBCI-1} (\forall \x, \y, \z\in \LomX)
   (\jI\jel (\x\op \y)\op ((\y\op \z)\op (\x\op \z))),
 \\&\label{OBCI-2} (\forall \x, \y\in \LomX)
   (\jI\jel \x\op ((\x\op \y)\op \y)),
 \\&\label{OBCI-3} (\forall \x\in \LomX)
   (\jI\jel \x\op \x),
 \\&\label{OBCI-4} (\forall \x, \y\in \LomX)
   (\jI\jel \x\op \y, \, \jI\jel \y\op \x ~\Rightarrow ~\x=\y),
 \\&\label{OBCI-5} (\forall \x, \y\in \LomX)
   (\x\jel \y ~\Leftrightarrow ~\jI\jel \x\op \y),
 \\&\label{OBCI-6} (\forall \x, \y\in \LomX)
   (\jI \jel \x, \, \x\jel \y ~\Rightarrow ~\jI \jel \y).
\end{align}
\end{definition}

\begin{proposition}[\cite{221005}]\label{prp31-221009}
	Every OBCI-algebra \Xo0 satisfies:
	\begin{align}
		&\label{(a1)} (\forall \x\in \LomX) (\jI \op \x = \x). 
		\\&\label{(a2)} 
		(\forall \x, \y, \z\in \LomX) (\z\op (\y\op \x)=\y\op (\z\op \x)).
 	\\&\label{(a3)} 
 	(\forall \x, \y, \z\in \LomX) (\jI \jel \x\op \y ~\Rightarrow ~\jI\jel (\y\op \z)\op (\x\op \z)).
 	\\&\label{(a4)} 
 	(\forall \x, \y, \z\in \LomX) (\jI\jel \x\op \y, \, \jI\jel \y\op \z ~\Rightarrow ~\jI\jel \x \op \z).
		\\&\label{(b4)} 
		(\forall \x, \y, \z\in \LomX) (\jI\jel (\y\op \z)\op ((\x\op \y)\op (\x\op \z))).
 	\\&\label{(b5)} 		(\forall \x, \y, \z\in \LomX) 
 	(\jI \jel \x\op \y ~\Rightarrow ~\jI\jel (\z\op \x)\op (\z\op \y)).
%%%%%%%%%%%%%%%%%%%%%%%%%%%%%%%%%%%%%%%%%%%%%%%%%%%%%%%%%%%%%%%%%%%%%%
	\end{align}
\end{proposition}

For future convenience, \Xo0 represents the OBCI-algebra unless otherwise specified.

\begin{definition}[\cite{221005}] \label{def32-221005}	
	A subset $\LomA$ of $\LomX$  is called 
	\begin{enumerate}
		\item[$\bullet$]  a \textcolor{red}{\it subalgebra}  of \Xo0 if it satisfies:
		\begin{align}\label{(7)} 
			(\forall \x, \y\in \LomX) (\x, \y\in \LomA ~\Rightarrow ~\x\op \y\in \LomA).
		\end{align}
		\item[$\bullet$]  an \textcolor{red}{\it ordered subalgebra}  of \Xo0 if it satisfies:
		\begin{align}\label{(77)}
			(\forall \x, \y\in \LomX) (\x, \y\in \LomA, \, \jI\jel \x, \, \jI \jel \y ~\Rightarrow ~\x\op \y\in \LomA).
		\end{align}
	\end{enumerate}	
\end{definition}

 \begin{definition}[\cite{221005}] \label{def33-221005}
 A subset $\LomF$ of $\LomX$ is called 
 \begin{enumerate}
 	\item[$\bullet$] a \textcolor{red}{\it filter}  of \Xo0 if it satisfies:
 	\begin{align}
 	&\label{flt-1} \jI \in \LomF,
 	\\&\label{flt-2}	(\forall \x, \y\in \LomX) 
    (\x\op \y \in \LomF, \, \x\in \LomF ~\Rightarrow ~\y\in \LomF).	
 	\end{align}
 	\item[$\bullet$] an \textcolor{red}{\it ordered filter} of \Xo0 if it satisfies \eqref{flt-1} and
 	\begin{align}	
 	&\label{Oflt} 
 		(\forall \x, \y\in \LomX)
 	(\x\in \LomF, \, \jI \jel \x\op \y ~\Rightarrow ~\y\in \LomF).
 	\end{align}
 \end{enumerate}
 \end{definition}

\begin{proposition}[\cite{221005}] \label{TYGV70}
			If an ordered filter $\LomF$ of \Xo0 satisfies 
 	\begin{align}	
 	\label{qcOT35-221215-1} 
 		(\forall \x\in \LomX)
 	(\x\in \LomF ~\Rightarrow ~\jI \jel \x),
 	\end{align}
then it is a filter of \Xo0.
\end{proposition} 

\begin{definition}[\cite{221005}] \label{Dq-closed}
	An (ordered) filter $\LomF$ of \Xo0 is said to be 
	\textcolor{red}{\it closed} if it is an (ordered) subalgebra of \Xo0.
\end{definition}

\begin{definition}[\cite{221005}] \label{def-dproduct}
Let $\textbf{\LomX} := (\LomX,$ $\opX,$ $\jIX,$ $\leX)$ and 
$\textbf{\LomY} := (\LomY,$ $\opY,$ $\jIY,$ $\leY)$
 be OBCI-algebras. Consider a binary operation $``\oq"$, 
a constant $``{\bf \jI}"$ and a binary relation $``\ll"$  in
the Cartesian product $\LomX \times \LomY$ defined as follows:
\begin{align*}
&(\x, \jx)\oq (\y, \jy)=(\x\opX \y, \jx\opY \jy),	
\\& {\bf \jI}=(\jIX, \jIY),
\\&  (\x, \jx)\ll (\y, \jy) ~\Leftrightarrow ~\x\leX \y, \, \jx \leY \jy)
\end{align*}
for all $(\x, \jx),  (\y, \jy)\in \LomX \times \LomY$. 
$\textbf{\LomX} \times \textbf{\LomY} := (\LomX \times \LomY,$ $\oq,$ ${\bf \jI},$ $\ll)$ is said to be a \textcolor{red}{\it direct product OBCI-algebra} if it
is an OBCI-algebra.
\end{definition}

%%%%%%%%%%%%%%%%%%%%%%%%%%%%%%%%%%%%%%%%%%%%%%%%%%%%%%%%%%%%%%%%%%%%%
\section{Ordered homomorphisms} 
%%%%%%%%%%%%%%%%%%%%%%%%%%%%%%%%%%%%%%%%%%%%%%%%%%%%%%%%%%%%%%%%%%%%%

In this section,  we first introduce ordered homomorphisms and kernels and then investigate several properties related to (ordered) subalgebras, (ordered) filters and direct products.

%%%%%%%%%%%%%%%%%%%%%%%%%%%%%%%%%%%%%%%%%%%%%%%%%%%%%%%%%%%%%%%%%%%%%%%%%%
\subsection{Ordered homomorphisms and kernels} %%%%%%%%%%%%%%%%%%%%%%%%%%%%%%%%%%%%%%%%%%%%%%%%%%%%%%%%%%%%%%%%%%%%%%%%%%

\begin{definition}\label{DqD161-230520}
	Let $\textbf{\LomX} := (\LomX,$ $\opX,$ $\jIX,$ $\leX)$ and
	$\textbf{\LomY} := (\LomY,$ $\opY,$ $\jIY,$ $\leY)$
	be OBCI-algebras. A mapping $\kap :\LomX \rightarrow \LomY$ is called 
 a 
\textcolor{red}{\it homomorphism} if it satisfies
\begin{align}\label{homo}
	(\forall \x, \y\in \LomX)(\kap(\x\opX \y)= \kap(\x)\opY \kap(\y));
\end{align}
	an \textcolor{red}{\it ordered map} (briefly, {\it $O$-map}) if it satisfies
	\begin{align}\label{O-map}
		(\forall \x, \y\in \LomX)(\jIX\leX \x\opX \y 
		~\Rightarrow ~\jIY \leY \kap(\x)\opY \kap(\y));
	\end{align}
	an \textcolor{red}{\it ordered homomorphism} (briefly, {\it $O$-homomorphism}) if it satisfies both \eqref{homo} and \eqref{O-map}.
\end{definition}

Given OBCI-algebras 
$\textbf{\LomX} := (\LomX,$ $\opX,$ $\jIX,$ $\leX)$ and
$\textbf{\LomY} := (\LomY,$ $\opY,$ $\jIY,$ $\leY)$, 
the mapping 
$\kap :\LomX \rightarrow \LomY, ~\x \mapsto \jIY,$
is an $O$-homomorphism. 

\begin{example}\label{EqE161-230520}
	Let $\LomX:=\{\jI, \x, \y\}$ and $\LomY=\{\jI,\jx\}$ be sets with  binary operations $``\opX"$ and $``\opY"$
given by Table \ref{[B3-2]} and Table \ref{[B2-1]}, respectively. 
%%%%%%%%%%%%%%%%%%%%%%%%%%%%%%%%%%%%%%%%%%%%%%%%%%%%%%%%%%%%%%%%%%%%%%%%%%
\begin{table}[!ht]
	\caption{Cayley table for the binary operation ``$\op$''}
	\vspace*{-8pt}
%%%%%%%%%%%%%%%%%%%%%%%%%%%%%%%%%%%%%%%%%%%%%%%%%%%%%%%%%%%%%%%%%%%%%%%%%%
	\begin{center}
		\def\temptablewidth{0.5\textwidth}
		{\rule{\temptablewidth}{1pt}}    \label{[B3-2]}
		\begin{tabular*}{\temptablewidth}{@{\extracolsep{\fill}}c|ccc}
			~$\opX$ &$\jI$ &$\x$ &$\y$~\\
			\hline
			~$\jI$ &$\jI$ &$\x$ &$\y$~\\
			~$\x$ &$\jI$ &$\jI$ &$\y$~\\
			~$\y$ &$\y$ &$\y$ &$\jI$~\\
		\end{tabular*}
		{\rule{\temptablewidth}{1pt}}
	\end{center}
%%%%%%%%%%%%%%%%%%%%%%%%%%%%%%%%%%%%%%%%%%%%%%%%%%%%%%%%%%%%%%%%%%%%%%%%%%
\end{table}
%%%%%%%%%%%%%%%%%%%%%%%%%%%%%%%%%%%%%%%%%%%%%%%%%%%%%%%%%%%%%%%%%%%%%%%%%%
\begin{table}[!ht]
	\caption{Cayley table for the binary operation ``$\op$''}
	\vspace*{-8pt}
%%%%%%%%%%%%%%%%%%%%%%%%%%%%%%%%%%%%%%%%%%%%%%%%%%%%%%%%%%%%%%%%%%%%%%%%%%
	\begin{center}
		\def\temptablewidth{0.5\textwidth}
		{\rule{\temptablewidth}{1pt}}    \label{[B2-1]}
		\begin{tabular*}{\temptablewidth}{@{\extracolsep{\fill}}c|cc}
			~$\opY$ &$\jI$ &$\jx$~\\
			\hline
			~$\jI$ &$\jI$ &$\jx$~\\
			~$\jx$ &$\jx$ &$\jI$~\\
		\end{tabular*}
		{\rule{\temptablewidth}{1pt}}
	\end{center}
%%%%%%%%%%%%%%%%%%%%%%%%%%%%%%%%%%%%%%%%%%%%%%%%%%%%%%%%%%%%%%%%%%%%%%%%%%
\end{table}
%%%%%%%%%%%%%%%%%%%%%%%%%%%%%%%%%%%%%%%%%%%%%%%%%%%%%%%%%%%%%%%%%%%%%%%%%%

\noindent
Let $\leX :=\{(\jI, \jI), (\x,\x), (\y,\y), (\x,\jI)\}$ and $\leY :=\{(\jI, \jI), (\jx, \jx)\}$. Then 
  $\textbf{\LomX} := (\LomX,$ $\opX,$ $\jI,$ $\leX)$ and
$\textbf{\LomY} := (\LomY,$ $\opY,$ $\jI,$ $\leY)$
are OBCI-algebras.  
Define a mapping 
\begin{align} \label{qcE161-221020}
\kap :\LomX \rightarrow \LomY,	~\z \mapsto \left\{\begin{array}{ll}
		\jI &\text{\rm if $\z\in \{\jI, \x\}$},\\
		\jx &\text{\rm if $\z =\y$}.
	\end{array}\right. 
\end{align}
It is routine to verify that $\kap$ is an $O$-homomorphism.
\end{example}

\begin{proposition}\label{LGV70-230520}
Let $\kap$ be a mapping from an OBCI-algebra $\textbf{\LomX} := (\LomX,$ $\opX,$ $\jIX,$ $\leX)$ to an OBCI-algebra
	$\textbf{\LomY} := (\LomY,$ $\opY,$ $\jIY,$ $\leY)$.
Let $\textbf{\LomX} := (\LomX,$ $\opX,$ $\jIX,$ $\leX)$ and
$\textbf{\LomY} := (\LomY,$ $\opY,$ $\jIY,$ $\leY)$
be OBCI-algebras.
If $\kap :\LomX \rightarrow \LomY$ is an $O$-homomorphism, then 
\begin{align}
		&\label{hqcP162-221014-1} \jIY\leY \kap(\jIX)\opY \kap(\jIX),
		\\&\label{hqcP162-221014-2} \jIY \leY \kap(\jIX).
		\\&\label{qcP162-221014-2} (\forall \x, \y\in \LomX)
			(\x\leX \y ~\Rightarrow \kap(\x)\leY \kap(\y)).
	\end{align}
\end{proposition}

\begin{proof}
\eqref{hqcP162-221014-1} and \eqref{hqcP162-221014-2} are straightforward by the combination of \eqref{OBCI-3},  \eqref{homo} and \eqref{(a1)}. \eqref{qcP162-221014-2} is straightforward by the combination of \eqref{OBCI-5} and \eqref{O-map}.
\end{proof}

The following examples show that homomorphisms and $O$-maps between OBCI-algebras may not be $O$-homomorphisms.

\begin{example}\label{ETkcGV70-230520}
Let $\LomX:=\{\1, \3, \0\}$ be a set with the binary operations $``\opX"$ given by Table \ref{[Knof]}. 
%%%%%%%%%%%%%%%%%%%%%%%%%%%%%%%%%%%%%%%%%%%%%%%%%%%%%%%%%%%%%%%%%%%%%%%%%%
\begin{table}[!ht]
	\caption{Cayley table for the binary operation ``$\op$''}
	\vspace*{-8pt}
%%%%%%%%%%%%%%%%%%%%%%%%%%%%%%%%%%%%%%%%%%%%%%%%%%%%%%%%%%%%%%%%%%%%%%%%%%
	\begin{center}
		\def\temptablewidth{0.5\textwidth}
		{\rule{\temptablewidth}{1pt}}    \label{[Knof]}
		\begin{tabular*}{\temptablewidth}{@{\extracolsep{\fill}}c|ccc}
			~$\opX$ &$\1$ &$\3$ &$\0$~\\
			\hline
			~$\1$ &$\1$ &$\0$ &$\0$~\\
			~$\3$ &$\1$ &$\3$ &$\0$~\\
			~$\0$ &$\1$ &$\1$ &$\1$~\\
		\end{tabular*}
		{\rule{\temptablewidth}{1pt}}
	\end{center}
%%%%%%%%%%%%%%%%%%%%%%%%%%%%%%%%%%%%%%%%%%%%%%%%%%%%%%%%%%%%%%%%%%%%%%%%%%
\end{table}
\noindent Let $\leX :=\{(\1, \1), (\3,\3), (\0,\0), (\3,\0), (\0, \3)\}$. Then 
  $\textbf{\LomX} := (\LomX,$ $\opX,$ $\3,$ $\leX)$ is an OBCI-algebra.  
Define a mapping 
\begin{align} \label{qcE161-230520}
\kap :\LomX \rightarrow \LomX,	~\z \mapsto \left\{\begin{array}{ll}
		\1 &\text{\rm if $\z= \0$},\\
		\3 &\text{\rm if $\z= \3$},\\
		\0 &\text{\rm if $\z =\1$}.
	\end{array}\right. 
\end{align}
It is routine to verify that $\kap$ is a homomorphism. However, $\kap$ is not an $O$-map since for $\1, \3 \in \LomX$, we have 
$\3 \leX \3 \opX \1$ but 
$$\3 = \kap(\3) \not\leX \0 = \3 \opX \0 = \kap(\3) \opY \kap(\1).$$
Therefore, the map $\kap$ does not form an $O$-homomorphism.
\end{example}

\begin{example}\label{EqvO2h-230429}
Let $\LomX:=\{\1, \jI, \p, \0\}$  and $\LomY:=\{\1, \8, \9, \0\}$ be sets with  binary operation $``\opX"$ 
and $``\opY"$ given by Table \ref{Edef31-1-221006-tbl} and Table \ref{[Y-230429]}, respectively. 
%%%%%%%%%%%%%%%%%%%%%%%%%%%%%%%%%%%%%%%%%%%%%%%%%%%%%%%%%%%%%%%%%%%%%%%%%%
%%%%%%%%%%%%%%%%%%%%%%%%%%%%%%%%%%%%%%%%%%%%%%%%%%%%%%%%%%%%%%%%%%%%%%%%%%
\begin{table}[!ht]
	\caption{Cayley table for the binary operation ``$\opX$''}
	\vspace*{-8pt}
%%%%%%%%%%%%%%%%%%%%%%%%%%%%%%%%%%%%%%%%%%%%%%%%%%%%%%%%%%%%%%%%%%%%%%%%%%
		\begin{center}
			\def\temptablewidth{0.5\textwidth}
			{\rule{\temptablewidth}{1pt}}    \label{Edef31-1-221006-tbl}
			\begin{tabular*}{\temptablewidth}{@{\extracolsep{\fill}}ccccc}
				~$\opX$ &$\1$ &$\jI$ &$\p$ &$\0$~\\
				\hline
				~$\1$ &$\1$ &$\0$ &$\0$ &$\0$~\\
				~$\jI$ &$\1$ &$\jI$ &$\p$ &$\0$~\\
				~$\p$ &$\1$ &$\p$ &$\jI$ &$\0$~\\
				~$\0$ &$\1$ &$\1$ &$\1$ &$\1$~\\
			\end{tabular*}
			{\rule{\temptablewidth}{1pt}}
		\end{center}
%%%%%%%%%%%%%%%%%%%%%%%%%%%%%%%%%%%%%%%%%%%%%%%%%%%%%%%%%%%%%%%%%%%%%%%%%%
\end{table}
%%%%%%%%%%%%%%%%%%%%%%%%%%%%%%%%%%%%%%%%%%%%%%%%%%%%%%%%%%%%%%%%%%%%%%%%%%
\begin{table}[!ht]
	\caption{Cayley table for the binary operation ``$\opY$''}
	\vspace*{-8pt}
%%%%%%%%%%%%%%%%%%%%%%%%%%%%%%%%%%%%%%%%%%%%%%%%%%%%%%%%%%%%%%%%%%%%%%%%%%
	\begin{center}
		\def\temptablewidth{0.5\textwidth}
		{\rule{\temptablewidth}{1pt}}    \label{[Y-230429]}
		\begin{tabular*}{\temptablewidth}{@{\extracolsep{\fill}}c|cccc}
		  ~$\opY$ &$\1$ &$\8$ &$\9$ &$\0$~\\
			\hline
			~$\1$ &$\1$ &$\0$ &$\0$   &$\0$~\\
			~$\8$ &$\1$ &$\8$ &$\9$   &$\0$~\\
			~$\9$ &$\1$ &$\8$   &$\8$ &$\0$~\\
			~$\0$ &$\1$   &$\1$   &$\1$   &$\1$~\\
		\end{tabular*}
		{\rule{\temptablewidth}{1pt}}
	\end{center}
%%%%%%%%%%%%%%%%%%%%%%%%%%%%%%%%%%%%%%%%%%%%%%%%%%%%%%%%%%%%%%%%%%%%%%%%%%
\end{table}

Let $$\leX :=\{(\0,\0), (\jI, \jI), (\p, \p), (\1, \1), (\0, \jI), (\0, \p), (\jI, \1), (\p, \1)\}$$ and
$$\leY :=\{(\1, \1), (\8, \8), (\9,\9), (\0,\0), (\8,\1), (\9, \8), (\0, \9)\}.$$ 
Then  $\textbf{\LomX} := (\LomX,$ $\opX,$ $\jI,$ $\leX)$ and
$\textbf{\LomY} := (\LomY,$ $\opY,$ $\8,$ $\leY)$ are OBCI-algebras. 
Define a mapping $\kap$ from $\LomX$ to $\LomY$ as follows:
\begin{align} \label{qcE54-230111}
	\kap :\LomX \rightarrow \LomY,	~w \mapsto \left\{\begin{array}{ll}
		\1 &\text{\rm if $w=\1$},\\
		\8 &\text{\rm if $w=\jI$},\\
    	\9 &\text{\rm if $w =\p$},\\
		\0 &\text{\rm if $w =\0$}.
	\end{array}\right. 
\end{align}
It is clear that $\kap$ satisfies \eqref{O-map}. However, as the following example shows, $\kap$ is not a homomorphism. For $\jI, \p \in \LomX$, we have 
$$\kap(\p \opX \jI) = \kap(\p) = \9 \neq \8 = \9 \opY \8 = \kap(\p) \opY \kap(\jI).$$ Therefore, the map $\kap$ does not form an $O$-homomorphism.
\end{example}

%%%%%%%%%%%%%%%%%%%%%%%%%%%%%%%%%%%%%%%%%%%%%%%%%%%%%%%%%%%%%%%%%%%%%%%%%%

\begin{definition}[\cite{221215}] \label{Dker-230108}
	Let $\kap$ be a mapping from an OBCI-algebra $\textbf{\LomX} := (\LomX,$ $\opX,$ $\jIX,$ $\leX)$ to an OBCI-algebra
	$\textbf{\LomY} := (\LomY,$ $\opY,$ $\jIY,$ $\leY)$.
The \textcolor{red}{\it kernel} of $\kap$ is defined to be a subset, denoted by ${\rm ker}(\kap)$, of $\LomX$ that satisfies:
\begin{align}\label{ker}	
	(\forall \x \in \LomX) (\x \in {\rm ker}(\kap) ~\Leftrightarrow ~ \jIY \leY \kap(\x)).	
\end{align}
\end{definition}

It is clear that ${\rm ker}(\kap):=\{\x\in \LomX \mid \jIY \leY \kap(\x)\}$, and it is unique. 

\begin{example}\label{EqEY-230520}
Let $\kap$ be a mapping in Example \ref{EqvO2h-230429}. 	As already mentioned in it, $\kap$ is not an $O$-homomorphism.
It is routine to calculate that ${\rm ker}(\kap)=\{\1, \jI\}$.
\end{example}

\begin{proposition}[\cite{221215}] \label{Pker-230106}
	If $\kap$ is a mapping from an OBCI-algebra $\textbf{\LomX} := (\LomX,$ $\opX,$ $\jIX,$ $\leX)$ to an OBCI-algebra
$\textbf{\LomY} := (\LomY,$ $\opY,$ $\jIY,$ $\leY)$, then
its kernel ${\rm ker}(\kap)$ is given by the following set.
	\begin{align}\label{ker*-230106}
		{\rm ker}(\kap)=\{\y\in \LomX \mid (\exists\x\in \LomX) 
		(\jIY \leY \kap(\x), \, \jIY \leY \kap(\x) \opY \kap(\y))\}.
	\end{align}
\end{proposition}

\begin{example}\label{(EX-OHK)} 
Let $\kap$ be an $O$-homomorphism in Example \ref{EqE161-230520}. We have 
 ${\rm ker}(\kap)= \{\jI\}$. 
\end{example}

%%%%%%%%%%%%%%%%%%%%%%%%%%%%%%%%%%%%%%%%%%%%%%%%%%%%%%%%%%%%%%%%%%%%%%%%%%
\subsection{Kernels and (ordered) subalgebras} %%%%%%%%%%%%%%%%%%%%%%%%%%%%%%%%%%%%%%%%%%%%%%%%%%%%%%%%%%%%%%%%%%%%%%%%%%

We first verify that kernels may not be (ordered) subalgebras.

\begin{example}\label{EGV80-230520}
Let $\textbf{\LomX} := (\LomX,$ $\opX,$ $\3,$ $\leX)$ be an OBCI-algebra given in Example \ref{ETkcGV70-230520}. 
Let $\kap :\LomX \rightarrow \LomX$ be an automorphism as the identity map. Then
 ${\rm ker}(\kap)=\{\1, \3\}$ and it is neither a subalgebra of $\textbf{\LomX}$  nor an ordered subalgebra of $\textbf{\LomX}$ since 
$\1 \op \3 = \0 \notin {\rm ker}(\kap)$.
\end{example}

\begin{definition}\label{DkOh-closed}
Let $\textbf{\LomX} := (\LomX,$ $\opX,$ $\jIX,$ $\leX)$ and
		$\textbf{\LomY} := (\LomY,$ $\opY,$ $\jIY,$ $\leY)$
		be OBCI-algebras, and $\kap :\LomX \rightarrow \LomY$ be an $O$-homomorphism. 
The kernel  of $\kap$ is said to be 
		\textcolor{red}{\it closed} if it is a subalgebra  of $\textbf{\LomX}$; and		\textcolor{red}{\it ordered closed} (briefly $O$-closed) if it is an ordered subalgebra  of $\textbf{\LomX}$.
\end{definition}

By `($O$-)closed,' we ambiguously denote the words `closed' and `$O$-closed' together if we need not distinguish them.

\begin{example}\label{Eq[E161]-230520} 
Let $\textbf{\LomX} := (\LomX,$ $\opX,$ $\jI,$ $\leX)$ be an OBCI-algebra in Example \ref{EqE161-230520}. Let $\kap :\LomX \rightarrow \LomX$ be an automorphism as the identity map. 
Then ${\rm ker}(\kap) = \{\jI\}$ and it is the ($O$-)closed kernel of $\kap$.
\end{example}

\begin{proposition}
\label{PGV60-230521}
	Let $\textbf{\LomX} := (\LomX,$ $\opX,$ $\jIX,$ $\leX)$ and
		$\textbf{\LomY} := (\LomY,$ $\opY,$ $\jIY,$ $\leY)$
		be OBCI-algebras, and $\kap :\LomX \rightarrow \LomY$ be an $O$-homomorphism. 
\begin{enumerate}
		\item[{\rm (i)}] If  ${\rm ker}(\kap)$  is the closed kernel of $\kap$, then it satisfies:
		\begin{align}
			\label{Tkcq-closed-J} 
			(\forall \x\in \LomX)
			(\kap(\jIX) \leY \kap(\x) ~\Rightarrow ~\x\opX \jIX \in {\rm ker}(\kap)).
		\end{align}
		\item[{\rm (ii)}] If ${\rm ker}(\kap)$  is the $O$-closed kernel of $\kap$ and satisfies \eqref{qcOT35-221215-1}, then it further satisfies \eqref{Tkcq-closed-J}.
 	\end{enumerate}

\end{proposition}

\begin{proof}
(i) Let ${\rm ker}(\kap)$ be the closed kernel of $\kap$ and $\x\in \LomX$ be such that $\kap(\jIX) \leY \kap(\x)$. Then $\jIY \leY \kap(\x)$ by \eqref{hqcP162-221014-2} and \eqref{OBCI-6}, and so $\x\in {\rm ker}(\kap)$ by \eqref{ker}. Since 
$\jIX\in {\rm ker}(\kap)$, we then have $\x\opX \jIX \in {\rm ker}(\kap)$	 by \eqref{(7)}.

(ii) Let ${\rm ker}(\kap)$ is the $O$-closed kernel of $\kap$ and $\x\in \LomX$ be such that $\kap(\jIX) \leY \kap(\x)$. Then, as in (i), $\jIY \leY \kap(\x)$ and so $\x\in {\rm ker}(\kap)$. Hence $\jIX \leX \x$ by \eqref{qcOT35-221215-1}. Then since $\jIX\in {\rm ker}(\kap)$ and $\jIX \leX \jIX$, we obtain $\x\opX \jIX \in {\rm ker}(\kap)$ by \eqref{(77)}.
\end{proof}

As seen in the following example, if the kernel ${\rm ker}(\kap)$ of $\kap$ is not ($O$-)closed in Proposition \ref{PGV60-230521}, then the condition \eqref{Tkcq-closed-J} does not hold.

\begin{example}\label{kEGV80-230117}
The kernel ${\rm ker}(\kap)=\{\1, \3\}$ in Example \ref{Eq[E161]-230520} does not satisfy \eqref{Tkcq-closed-J} 
since $\kap(\3) \jel \kap(\1)$ but $\1 \op \3 = \0 \notin {\rm ker}(\kap)$.
\end{example}

The following theorem verifies that the assertion \eqref{Tkcq-closed-J} is a condition for the kernel to be ($O$-)closed.

\begin{theorem}\label{Tkq-(o)closed}
Let $\textbf{\LomX} := (\LomX,$ $\opX,$ $\jIX,$ $\leX)$ and $\textbf{\LomY} := (\LomY,$ $\opY,$ $\jIY,$ $\leY)$ be OBCI-algebras, and $\kap :\LomX \rightarrow \LomY$ be an $O$-homomorphism satisfying $\jIY = \kap(\jIX)$. 
If the kernel ${\rm ker}(\kap)$ of $\kap$ satisfies \eqref{Tkcq-closed-J}, 
		then it is ($O$-)closed.
\end{theorem}

\begin{proof} 
	Let ${\rm ker}(\kap)$ satisfy \eqref{Tkcq-closed-J}. First it is proved that ${\rm ker}(\kap)$ is closed. We verify that ${\rm ker}(\kap)$ is a subalgebra of $\textbf{\LomX}$. Let $\x\in {\rm ker}(\kap)$ and $\y\in {\rm ker}(\kap)$. 
Then $\x \op \jI \in {\rm ker}(\kap)$ and $\y \op \jI \in {\rm ker}(\kap)$ by \eqref{Tkcq-closed-J}. Thus $\jIY \leY \kap(\x)$, $\jIY \leY \kap(\y)$, $\jIY \leY \kap(\x) \opY \jIY$ and $\jIY \leY \kap(\y) \opY \jIY$ by \eqref{ker} and $\jIY = \kap(\jIX)$. Moreover $\jIX \opX \x=\x\in {\rm ker}(\kap)$ and $\jIX \opX \y=\y \in {\rm ker}(\kap)$  by \eqref{(a1)}, and so $\jIY \leY \jIY \opY \kap(\x) = \kap(\x)$ and $\jIY \leY \jIY \opY \kap(\y) = \kap(\y)$ by \eqref{ker} and $\jIY = \kap(\jIX)$. 
		Note that 
$$\jIY \leY (\kap(\x) \opY \jIY) \opY ((\jIY \opY \kap(\y)) \opY (\kap(\x) \opY \kap(\y)))$$
		by \eqref{OBCI-1}. Then it follows from \eqref{OBCI-6} that 
$$\jIY \leY (\jIY \opY \kap(\y)) \opY (\kap(\x) \opY \kap(\y))$$ and so 
$$\jIY \leY \kap(\y) \opY (\kap(\x) \opY \kap(\y)).$$ 
Analogously we further obtain $\jIY \leY \kap(\x) \opY \kap(\y)$, and so 
$$\jIY \leY \kap(\x) \opY \kap(\y) = \kap(\x \opX \y)$$ 
by \eqref{homo}. Hence $\x \op \y \in {\rm ker}(\kap)$ by \eqref{ker}, and therefore ${\rm ker}(\kap)$ is a subalgebra of $\textbf{\LomX}$.

The proof for the case that ${\rm ker}(\kap)$ is $O$-closed is almost the same..\end{proof}

Let the kernel ${\rm ker}(\kap)$ of $\kap$ be not ($O$-)closed in Proposition \ref{PGV60-230521}. Then the condition \eqref{Tkcq-closed-J} does not hold as the following example shows. 

\begin{example}\label{kEGV80-230520}
The kernel ${\rm ker}(\kap)=\{\1, \3\}$ in Example \ref{EGV80-230520} does not satisfy \eqref{Tkcq-closed-J} 
since $\3 = \kap(\3) \jel \kap(\1) = \1$ but $\1 \op \3 = \0 \notin {\rm ker}(\kap)$.
\end{example}

%%%%%%%%%%%%%%%%%%%%%%%%%%%%%%%%%%%%%%%%%%%%%%%%%%%%%%%%%%%%%%%%%%%%%%

\begin{theorem}\label{TohGV50-230520}
	Let $\kap :\LomX \rightarrow \LomY$ be an $O$-homomorphism from an OBCI-algebra  
		$\textbf{\LomX} := (\LomX,$ $\opX,$ $\jIX,$ $\leX)$ to an OBCI-algebra 
		$\textbf{\LomY} := (\LomY,$ $\opY,$ $\jIY,$ $\leY)$.
		\begin{enumerate}
			\item[\rm (i)] If $\LomG$ is a subalgebra of $\textbf{\LomY}$,  
			then $\kap^{-1}(\LomG)$ is a subalgebra of $\textbf{\LomX}$.
\item[\rm (ii)] Suppose that $\kap$ is surjctive.
			If $\LomF$ is a subalgebra of $\textbf{\LomX}$,  
			then $\kap(\LomF)$ 	is a subalgebra of $\textbf{\LomY}$.
		\end{enumerate}
\end{theorem}

\begin{proof}
(i) Assume that $\x, \y\in \kap^{-1}(\LomG)$. We have $\kap(\x)\in \LomG$ and $\kap(\y)\in \LomG$, and thus $\kap(\x)\opY \kap(\y) \in \LomG$ by \eqref{(7)}. Then $\kap(\x\opX \y) \in \LomG$ by \eqref{homo}, and so $\x\opX \y\in \kap^{-1}(\LomG)$. Therefore $\kap^{-1}(\LomG)$ is a subalgebra of $\textbf{\LomX}$.
	
(ii) Assume that $\kap$ is surjctive, $\jx \in \kap(\LomF)$ and $\jy \in \kap(\LomF)$. Since $\kap$ is surjctive, there are $\x, \y\in \LomF$ such that $\jx=\kap(\x)$ and $\jy=\kap(\y)$. Then $\x \opX \y\in \LomF$ by \eqref{(7)}, and so $\kap(\x\opX \y)\in \kap(\LomF)$. Hence $\kap(\x\opX \y) = \kap(\x) \opY \kap(\y) =  \jx \opY \jy \in \kap(\LomF)$ by \eqref{homo}. Therefore $\kap(\LomF)$ 	is a subalgebra of $\textbf{\LomY}$.
\end{proof}

\begin{theorem}\label{Toh2GV50-230520}
Let $\kap :\LomX \rightarrow \LomY$ be an $O$-homomorphism from an OBCI-algebra  
		$\textbf{\LomX} := (\LomX,$ $\opX,$ $\jIX,$ $\leX)$ to an OBCI-algebra 
		$\textbf{\LomY} := (\LomY,$ $\opY,$ $\jIY,$ $\leY)$.
		\begin{enumerate}
			\item[\rm (i)] If $\LomG$ is an ordered subalgebra of $\textbf{\LomY}$, then $\kap^{-1}(\LomG)$ is an ordered subalgebra of $\textbf{\LomX}$.
			\item[\rm (ii)] Suppose that $\kap$ is surjective. 
			If $\LomF$ is an ordered subalgebra of $\textbf{\LomX}$ and satisfies \eqref{qcOT35-221215-1}, 
			then $\kap(\LomF)$ is an ordered subalgebra of $\textbf{\LomY}$.
		\end{enumerate}
\end{theorem}

\begin{proof}
(i) Let $\LomG$ be an ordered subalgebra of $\textbf{\LomY}$ and $\x, \y\in \kap^{-1}(\LomG)$ be such that $\jIX \leX \x$ and $\jIX \leX \y$. Then $\kap(\x)\in \LomG$ and $\kap(\y)\in \LomG$, and $$\jIY \leY \kap(\jIX) \leY \kap(\x)$$ and $$\jIY \leY \kap(\jIX) \leY \kap(\y)$$ 
by \eqref{qcP162-221014-2} and \eqref{hqcP162-221014-2}. Thus $\jIY \leY \kap(\x)$ and $\jIY \leY \kap(\y)$, and so $\kap(\x)\opY \kap(\y) \in \LomG$ by \eqref{(77)}. Hence
$$\kap(\x\opX \y)=\kap(\x)\opY \kap(\y) \in \LomG$$
by \eqref{homo}, and so $\x\opX \y\in \kap^{-1}(\LomG)$. Therefore $\kap^{-1}(\LomG)$ is an ordered subalgebra of $\textbf{\LomX}$.
	
(ii) Let $\kap$ be surjective, and $\LomF$ be an ordered subalgebra of $\textbf{\LomX}$ and satisfy \eqref{qcOT35-221215-1}. Let $\jx, \jy \in \LomY$ be such that $\jIY \leY \jx \in \kap(\LomF)$ and $\jIY \leY \jy \in \kap(\LomF)$. Then there are $\x, \y\in \LomF$ such tha $\kap(\x)=\jx$ and $\kap(\y)=\jy$ by the surjectivity of $\kap$. Hence $\jIX \leX \x$ and $\jIX \leX \y$ by \eqref{qcOT35-221215-1}, and so $\x \opX \y \in \LomF$ by \eqref{(77)}. Hence $\kap(\x\opX \y)\in \kap(\LomF)$, and so
$$\kap(\x)\opY \kap(\y)=\jx \opY \jy \in \kap(\LomF)$$ 
by \eqref{homo}. Therefore
		$\kap(\LomF)$ 	is an ordered subalgebra of $\textbf{\LomY}$.
\end{proof}

In Theorem \ref{Toh2GV50-230520}(ii), we may deal with such $O$-homomorphism between ordered subalgebras by adding some condition to the $O$-homomorphism of $\kap$ itself in place of adding \eqref{qcOT35-221215-1} to the ordered subalgebra $\LomF$ of $\textbf{\LomX}$.

\begin{theorem}\label{T2ohGV50-230523}
Let $\kap :\LomX \rightarrow \LomY$ be a surjective $O$-homomorphism from an OBCI-algebra  
		$\textbf{\LomX} := (\LomX,$ $\opX,$ $\jIX,$ $\leX)$ to an OBCI-algebra 
		$\textbf{\LomY} := (\LomY,$ $\opY,$ $\jIY,$ $\leY)$, and satisfies
		\begin{align}\label{Tcoh-230415}
			(\forall \x\in \LomX) (\jIY\leY \kap(\x) ~\Rightarrow ~\jIX \leX \x).
		\end{align} 
			If $\LomF$ is an ordered subalgebra of $\textbf{\LomX}$, then $\kap(\LomF)$ is an ordered subalgebra of $\textbf{\LomY}$.
\end{theorem}

\begin{proof}
Suppose that $\kap$ is surjective and satisfy \eqref{Tcoh-230415}. Let $\jx, \jy \in \LomY$ be such that $\jIY \leY \jx \in \kap(\LomF)$ and $\jIY \leY \jy \in \kap(\LomF)$. There are $\x, \y\in \LomF$ such that $\kap(\x)=\jx$ and $\kap(\y)=\jy$ by the surjectivity of $\kap$. Then $\jIX \leX \x$ and $\jIX \leX \y$ by \eqref{Tcoh-230415}, and so $\x \opX \y \in \LomF$ by \eqref{(77)}. Hence, as in the proof of Theorem \ref{Toh2GV50-230520}(ii), we obtain $\jx \opY \jy =\kap(\x)\opY \kap(\y) \in \kap(\LomF)$. Therefore 
		$\kap(\LomF)$ 	is an ordered subalgebra of $\textbf{\LomY}$.
\end{proof}

%%%%%%%%%%%%%%%%%%%%%%%%%%%%%%%%%%%%%%%%%%%%%%%%%%%%%%%%%%%%%%%%%%%%%%%%%%
\subsection{Kernels and (ordered) filters} %%%%%%%%%%%%%%%%%%%%%%%%%%%%%%%%%%%%%%%%%%%%%%%%%%%%%%%%%%%%%%%%%%%%%%%%%%

Let $\textbf{\LomX} := (\LomX,$ $\opX,$ $\jIX,$ $\leX)$ and $\textbf{\LomY} := (\LomY,$ $\opY,$ $\jIY,$ $\leY)$ be OBCI-algebras. If $\kap :\LomX \rightarrow \LomY$ is an $O$-homomorphism from 	$\textbf{\LomX}$ to $\textbf{\LomY}$, then the kernel of $\kap$ is an (ordered) filter of $\textbf{\LomX}$. 

\begin{theorem}\label{TkcGV70}
Let $\kap :\LomX \rightarrow \LomY$ be an $O$-homomorphism from an OBCI-algebra  
		$\textbf{\LomX} := (\LomX,$ $\opX,$ $\jIX,$ $\leX)$ to an OBCI-algebra 
		$\textbf{\LomY} := (\LomY,$ $\opY,$ $\jIY,$ $\leY)$. ${\rm ker}(\kap)$ is a filter of $\textbf{\LomX}$.
\end{theorem}

\begin{proof}
First we have $\jIX \in {\rm ker}(\kap)$ by \eqref{hqcP162-221014-2} and \eqref{ker}. Let $\x, \y\in \LomX$ be such that $\x \in {\rm ker}(\kap)$ and $\x \opX \y \in {\rm ker}(\kap)$. Then 
		$$\jIY \leY \kap(\x\opX \y)=\kap(\x)\opY \kap(\y)$$
by \eqref{ker} and \eqref{homo}, and so $\kap(\x)\leY \kap(\y)$ by \eqref{OBCI-5}. 
		Since $\jIY \leY \kap(\x)$ by \eqref{ker}, we obtain $\jIY \leY \kap(\y)$ by \eqref{OBCI-6}, and so $\y\in {\rm ker}(\kap)$.
		Therefore ${\rm ker}(\kap)$ is a filter of $\textbf{\LomX}$.
\end{proof}

\begin{theorem}\label{TOkcGV70}
If $\kap :\LomX \rightarrow \LomY$ is an $O$-homomorphism from an OBCI-algebra  
		$\textbf{\LomX} := (\LomX,$ $\opX,$ $\jIX,$ $\leX)$ to an OBCI-algebra 
		$\textbf{\LomY} := (\LomY,$ $\opY,$ $\jIY,$ $\leY)$, then ${\rm ker}(\kap)$ is an ordered filter of $\textbf{\LomX}$.		
\end{theorem}

\begin{proof}
It is clear that $\jIX \in {\rm ker}(\kap)$. 
Let $\x, \y\in \LomX$ be such that 	$\x \in {\rm ker}(\kap)$ and $\jIX \leX \x \opX \y$.
Then $\jIY \leY \kap(\x)$ by \eqref{ker} and $\jIY \leY \kap(\x) \opY \kap(\y)$ by \eqref{O-map}.
It follows from \eqref{OBCI-5} and \eqref{OBCI-6} that $\jIY \leY \kap(\y)$, that is, $\y \in {\rm ker}(\kap)$.
Hence ${\rm ker}(\kap)$ is an ordered filter of $\textbf{\LomX}$.
\end{proof}

%%%%%%%%%%%%%%%%%%%%%%%%%%%%%%%%%%%%%%%%%%%%%%%%%%%%%%%%%%%%%%%%%%%%%%%
We then display $O$-homomorphisms between between (ordered) filters of OBCI-algebras.

\begin{theorem}\label{TohqP165-230523}
Let $\kap :\LomX \rightarrow \LomY$ be an $O$-homomorphism from an OBCI-algebra  
		$\textbf{\LomX} := (\LomX,$ $\opX,$ $\jIX,$ $\leX)$ to an OBCI-algebra 
		$\textbf{\LomY} := (\LomY,$ $\opY,$ $\jIY,$ $\leY)$ such that $\jIY = \kap(\jIX)$.
		\begin{enumerate}
			\item[\rm (i)] If $\LomG$ is a  filter of $\textbf{\LomY}$, 
then $\kap^{-1}(\LomG)$
			is a  filter of $\textbf{\LomX}$.
			\item[\rm (ii)] Suppose that $\kap$ is surjective.
			If $\LomF$ is a  filter of $\textbf{\LomX}$,  then $\kap(\LomF)$ 	is a  filter of $\textbf{\LomY}$.
		\end{enumerate}
\end{theorem}

\begin{proof}
(i) Suppose that $\LomG$ is a  filter of $\textbf{\LomY}$. It is clear that $\jIX \in \kap^{-1}(\LomG)$ since $\jIY \in \LomG$ and $\jIY = \kap(\jIX)$. Let $\x\in \kap^{-1}(\LomG)$ and $\x\opX \y \in \kap^{-1}(\LomG)$ for $\x, \y\in \LomX$. Then $\kap(\x)\in \LomG$ and $\kap(\x)\opY \kap(\y)=\kap(\x\opX \y)\in \LomG$ by \eqref{homo}. Hence $\kap(\y)\in \LomG$ by \eqref{flt-2}, and so $\y \in \kap^{-1}(\LomG)$. Therefore 
		$\kap^{-1}(\LomG)$
		is a filter of $\textbf{\LomX}$.
	
(ii) Suppose that $\LomF$ is a  filter of $\textbf{\LomX}$. Then $\jIY \in \kap(\LomF)$ since $\jIX \in \LomF$ and $\jIY = \kap(\jIX)$. Let $\jx \in \kap(\LomF)$ and $\jx \opY \jy \in \kap(\LomF)$ for $\jx, \jy \in \LomY$. 
Then there are $\x\in \LomF$ and $\y\in \LomX$ such that $\kap(\x)=\jx$ and $\kap(\y)=\jy$ by the surjectivity of $\kap$. Hence 
$$\kap(\x\opX \y)=\kap(\x)\opY \kap(\y)=\jx \opY \jy \in \kap(\LomF)$$
by \eqref{homo}, and so $\x\opX \y \in \LomF$. Then $\y\in \LomF$ by \eqref{flt-2} and so $\jy =\kap(\y)\in \kap(\LomF)$. 
Therefore 
		$\kap(\LomF)$ 	is a filter of $\textbf{\LomY}$.
\end{proof}

The combination of Theorems \ref{TohGV50-230520} and \ref{TohqP165-230523} assures that 
if $\LomG$ is a  closed filter of $\textbf{\LomY}$,  then $\kap^{-1}(\LomG)$
is a  closed filter of $\textbf{\LomX}$; and 
if $\LomF$ is a closed filter of $\textbf{\LomX}$,  then $\kap(\LomF)$ 	is a closed filter of $\textbf{\LomY}$.

%%%%%%%%%%%%%%%%%%%%%%%%%%%%%%%%%%%%%%%%%%%%%%%%%%%%%%%%%%%%%%%%%%%%%%%

\begin{theorem}\label{T1of-230523}
Let $\kap :\LomX \rightarrow \LomY$ be an $O$-homomorphism from an OBCI-algebra  
		$\textbf{\LomX} := (\LomX,$ $\opX,$ $\jIX,$ $\leX)$ to an OBCI-algebra 
		$\textbf{\LomY} := (\LomY,$ $\opY,$ $\jIY,$ $\leY)$ such that $\jIY = \kap(\jIX)$.
		\begin{enumerate}
			\item[\rm (i)] If $\LomG$ is an ordered filter of $\textbf{\LomY}$, then $\kap^{-1}(\LomG)$ is an ordered filter of $\textbf{\LomX}$.
			\item[\rm (ii)] Suppose that $\kap$ is surjective  and satisfies \eqref{Tcoh-230415}.
			If $\LomF$ is an ordered filter of $\textbf{\LomX}$,  then $\kap(\LomF)$ 	is an ordered filter of $\textbf{\LomY}$.
		\end{enumerate}
\end{theorem}

\begin{proof}
(i) Suppose that $\LomG$ is an ordered filter of $\textbf{\LomY}$. Since $\jIY \in \LomG$ and $\jIY = \kap(\jIX)$, we have $\jIX \in \kap^{-1}(\LomG)$ by \eqref{flt-1}. Let $\x\in \kap^{-1}(\LomG)$ and $\jIX \leX \x\opX \y$ for $\x, \y\in \LomX$. Then $\kap(\x)\in \LomG$, and 
$$\jIY \leY \kap(\x)\opY \kap(\y)$$ 
by \eqref{O-map}. Hence $\kap(\y)\in \LomG$ by \eqref{Oflt}, and so $\y \in \kap^{-1}(\LomG)$. Therefore 
		$\kap^{-1}(\LomG)$
		is an ordered filter of $\textbf{\LomX}$.
	
(ii) Suppose that $\kap$ is surjective  and satisfies \eqref{Tcoh-230415}, and $\LomF$ is an ordered filter of $\textbf{\LomX}$. It is clear that $\jIY \in \kap(\LomF)$ as above.  Let $\jx \in \kap(\LomF)$ and $\jIY \leY \jx \opY \jy$ for $\jx, \jy \in \LomY$. Then there are $\x\in \LomF$ and $\y\in \LomX$ such that $\kap(\x)=\jx$ and $\kap(\y)=\jy$ by the surjectivity of $\kap$. Hence
$$\jIY \leY \jx \opY \jy = \kap(\x)\opY \kap(\y) = \kap(\x\opX \y)$$
by \eqref{homo}, and so $\jIX \leX \x\opX \y$ by \eqref{Tcoh-230415}. Then $\y\in \LomF$ by \eqref{Oflt}, and so $\jy =\kap(\y)\in \kap(\LomF)$. 
Therefore 
		$\kap(\LomF)$ 	is an ordered filter of $\textbf{\LomY}$.
\end{proof}

Let us take $\LomF$ as a filter containing the kernel of $\kap$ in Theorem \ref{T1of-230523}(ii). Then we can address the homomorphism between ordered filters as follows.

\begin{theorem}\label{T2of-230523}
Let $\kap :\LomX \rightarrow \LomY$ be an $O$-homomorphism from an OBCI-algebra  
	$\textbf{\LomX} := (\LomX,$ $\opX,$ $\jIX,$ $\leX)$ to an OBCI-algebra 
	$\textbf{\LomY} := (\LomY,$ $\opY,$ $\jIY,$ $\leY)$ such that $\jIY = \kap(\jIX)$.
Suppose that $\kap$ is surjective. If $\LomF$ is an ordered filter of $\textbf{\LomX}$, contains ${\rm ker}(\kap)$ and satisfies  \eqref{qcOT35-221215-1},
then 
		$\kap(\LomF)$ is an ordered  filter of $\textbf{\LomY}$.
\end{theorem} 

\begin{proof}
Suppose that $\kap$ is surjective, and $\LomF \subseteq \LomX$ contains ${\rm ker}(\kap)$ and satisfies \eqref{qcOT35-221215-1}. As above, $\jIY \in \kap(\LomF)$. 
	Let $\jx \in \kap(\LomF)$ and $\jIY \leY \jx \opY \jy$ for $\jx, \jy \in \LomY$. As above, there are $\x\in \LomF$ and $\y\in \LomX$ such that 
	$\kap(\x)=\jx \in \kap(\LomF)$ and $\kap(\y)=\jy$, and so $\jIY \leY \kap(\x) \opY \kap(\y)$. Moreover $\jIX \leX \jIX\opX \x$ by \eqref{qcOT35-221215-1} and \eqref{(a1)}, and so
$$\jIY \leY \kap(\jIX) \opY \kap(\x) = \kap(\x)$$
by \eqref{O-map}, $\jIY = \kap(\jIX)$ and \eqref{(a1)}.
Then  $\jIY \leY \kap(\y)$ by \eqref{OBCI-5} and \eqref{OBCI-6}, and so $\y\in {\rm ker}(\kap) \subseteq \LomF$. Hence $\jy =\kap(\y)\in \kap(\LomF)$. Therefore
	$\kap(\LomF)$ 	is an odered filter of $\textbf{\LomY}$.
%%%%%%%%%%%%%%%%%%%%%%%%%%%%%%%%%%%%%%%%%%%%%%%%%%%%%%%%%%%%%%%%%%	
\end{proof}

%%%%%%%%%%%%%%%%%%%%%%%%%%%%%%%%%%%%%%%%%%%%%%%%%%%%%%%%%%%%%%%%%%%%%%%

\begin{theorem}\label{TqT312-221114}
Let $\kap :\LomX \rightarrow \LomY$ be an $O$-homomorphism from an OBCI-algebra  
	$\textbf{\LomX} := (\LomX,$ $\opX,$ $\jIX,$ $\leX)$ to an OBCI-algebra 
	$\textbf{\LomY} := (\LomY,$ $\opY,$ $\jIY,$ $\leY)$ such that 
	$\kap$ is surjective and satisfies  $\jIY = \kap(\jIX)$.
	If we consider two sets:
	\begin{align*}
		&{\mathcal F}:=\{\LomF \mid \text{\rm $\LomF$ is a filter of $\LomX$ containing ${\rm ker}(\kap)$}\},
		\\&{\mathcal G}:=\{\LomG \mid \text{\rm $\LomG$ is a filter of $\LomY$}\},
	\end{align*}
	then there exists a bijective function 
	\begin{align}\label{qcT312-221114}
		\kbp : {\mathcal F} \rightarrow {\mathcal G}, ~\LomF \mapsto \kbp(\LomF)
	\end{align}
	such that $\kbp^{-1}(\LomG)=\kap^{-1}(\LomG)$. 	
\end{theorem}

\begin{proof}
We first show that 
	\begin{align}\label{qcT312-221114-1}
		(\forall \LomF \in {\mathcal F})(\kap^{-1}(\kap(\LomF))=\LomF).
	\end{align}
	Let $\LomF \in {\mathcal F}$. It is clear that $\LomF \subseteq \kap^{-1}(\kap(\LomF))$. 
	Let $\x\in \kap^{-1}(\kap(\LomF))$. Then $\kap(\x)\in \kap(\LomF)$
	and so $\kap(\x)=\kap(\z)$ for some $\z\in \LomF$. Then, since $\jIX \leX \x \opX \x$ by \eqref{OBCI-3}, we obtain
$$\jIY \leY \kap(\x) \opY \kap(\x) = \kap(\z) \opY \kap(\x) =  \kap(\z \opX \x)$$ 
by \eqref{O-map} and \eqref{homo}. 
Hence $\z \opX \x \in {\rm ker}(\kap)$ by \eqref{ker}. Thus we get $\z \opX \x \in \LomF$ since ${\rm ker}(\kap) \subseteq \LomF$. Then $\x \in \LomF$ by \eqref{flt-2}, and so $\kap^{-1}(\kap(\LomF))\subseteq \LomF$. Therefore \eqref{qcT312-221114-1} is valid. 
	Using Theorem \ref{TohqP165-230523}(ii), 
we can consider the mapping 
	$\kbp : {\mathcal F} \rightarrow {\mathcal G}, ~\LomF \mapsto \kbp(\LomF)$
	given by $\kbp(\LomF)=\kap(\LomF)$ for all $\LomF \in {\mathcal F}$. 
	Suppose that $\kbp(\LomF_1)=\kbp(\LomF_2)$ for all $\LomF_1, \LomF_2 \in {\mathcal F}$.
	Then $\kap(\LomF_1)=\kap(\LomF_2)$, which implies from \eqref{qcT312-221114-1} that
	$$\LomF_1=\kap^{-1}(\kap(\LomF_1))=\kap^{-1}(\kap(\LomF_2))=\LomF_2.$$
	Hence $\kbp$ is one-to-one. Since $\kap$ is a surjective $O$-homomorphism, we know that $\LomG =\kap(\kap^{-1}(\LomG))$ 
	and $\kap^{-1}(\LomG)$ is a filter of $\textbf{\LomX}$ containing ${\rm ker}(\kap)$ for all $\LomG \in {\mathcal G}$, i.e., $\kap^{-1}(\LomG) \in {\mathcal F}$,
	and $\kbp(\kap^{-1}(\LomG))=\LomG$. Thus $\kbp$ is onto. Therefore $\kbp$ is a bijective function and 
	$\kbp^{-1}(\LomG)=\kap^{-1}(\LomG)$. 
\end{proof}

\begin{theorem}\label{ToqT312-230501}
Let $\kap :\LomX \rightarrow \LomY$ be an $O$-homomorphism from an OBCI-algebra  
		$\textbf{\LomX} := (\LomX,$ $\opX,$ $\jIX,$ $\leX)$ to an OBCI-algebra 
		$\textbf{\LomY} := (\LomY,$ $\opY,$ $\jIY,$ $\leY)$ such that 
		$\kap$ is surjective and satisfies  $\jIY = \kap(\jIX)$.
		If we consider two sets:
		\begin{align*}
			&{\mathcal F}:=\{\LomF \mid \text{\rm $\LomF$ is an ordered filter of $\LomX$  containing ${\rm ker}(\kap)$ and satisfying \eqref{qcOT35-221215-1}}\},
			\\&{\mathcal G}:=\{\LomG \mid \text{\rm $\LomG$ is an ordered filter of $\LomY$}\},
		\end{align*}
		then there exists a bijective function $\kbp$ satisfying \eqref{qcT312-221114} as in Theorem \ref{TqT312-221114}.	
\end{theorem}

\begin{proof}
	Suppose that $\LomF \in {\mathcal F}$ contains ${\rm ker}(\kap)$ and satisfies \eqref{qcOT35-221215-1}. We first verify \eqref{qcT312-221114-1}. Clearly $\LomF \subseteq \kap^{-1}(\kap(\LomF))$. 
	Let $\x\in \kap^{-1}(\kap(\LomF))$. As in the proof of (ii) in Theorem \ref{TqT312-221114}, 
$$\jIY \leY \kap(\x) \opY \kap(\x) = \kap(\z) \opY \kap(\x).$$ 
Then 
$$\jIY \leY \kap(\z \opX \x) = \kap(\z) \opY \kap(\x)$$ 
by \eqref{homo}, and so $\z\opX \x \in {\rm ker}(\kap)$ by \eqref{ker}. Hence $\z\opX \x \in \LomF$, and so $\jIX \leX \z\opX \x$ by \eqref{qcOT35-221215-1}. Then $\x\in \LomF$ by \eqref{Oflt}, and so $\kap^{-1}(\kap(\LomF))\subseteq \LomF$. Therefore \eqref{qcT312-221114-1} is valid. 
	Using Theorem \ref{T2of-230523}, we can consider the mapping 
	$\kbp : {\mathcal F} \rightarrow {\mathcal G}, ~\LomF \mapsto \kbp(\LomF)$
	given by $\kbp(\LomF)=\kap(\LomF)$ for all $\LomF \in {\mathcal F}$. 
	Suppose that $\kbp(\LomF_1)=\kbp(\LomF_2)$ for all $\LomF_1, \LomF_2 \in {\mathcal F}$.
	Then $\kap(\LomF_1)=\kap(\LomF_2)$, which implies from \eqref{qcT312-221114-1} that
	$$\LomF_1=\kap^{-1}(\kap(\LomF_1))=\kap^{-1}(\kap(\LomF_2))=\LomF_2.$$
	Hence $\kbp$ is one-to-one. Since $\kap$ is a surjective homomorphism, we know that $\LomG =\kap(\kap^{-1}(\LomG))$ 
	and $\kap^{-1}(\LomG)$ is an ordered filter of $\textbf{\LomX}$
for all $\LomG \in {\mathcal G}$, i.e., $\kap^{-1}(\LomG) \in {\mathcal F}$,
	and $\kbp(\kap^{-1}(\LomG))=\LomG$. Thus $\kbp$ is onto. Therefore $\kbp$ is a bijective function and 
	$\kbp^{-1}(\LomG)=\kap^{-1}(\LomG)$. 
\end{proof}

%%%%%%%%%%%%%%%%%%%%%%%%%%%%%%%%%%%%%%%%%%%%%%%%%%%%%%%%%%%%%%%%%%%%%%%%%%
\subsection{Kernels and direct products} %%%%%%%%%%%%%%%%%%%%%%%%%%%%%%%%%%%%%%%%%%%%%%%%%%%%%%%%%%%%%%%%%%%%%%%%%%

Henceforth, we suppose that $\textbf{\LomX}_{1}$ := $(\LomX_{1}, 
\opA$, $\jIA$, $\leA)$, $\textbf{\LomX}_{2}$ := $(\LomX_{2}, \opB$, $\jIB$, $\leB)$, $\textbf{\LomY}_{1}$ := $(\LomY_{1}, \opC$, $\jIC$, $\leC)$ and $\textbf{\LomY}_{2}$ := $(\LomY_{2}, \opD$, $\jID$, $\leD)$ are OBCI-algebras.

\begin{definition}\label{def-ohdproduct}
Suppose that $\kap_{1} :\LomX_{1} \rightarrow \LomY_{1}$ and $\kap_{2} :\LomX_{2} \rightarrow \LomY_{2}$ are mappings from $\textbf{\LomX}_{1}$ to $\textbf{\LomY}_{1}$ and from $\textbf{\LomX}_{2}$ to 	$\textbf{\LomY}_{2}$, respectively. Consider the following direct product OBCI-algebras: 
\begin{center}
$\textbf{\LomX} = \textbf{\LomX}_{1} \times \textbf{\LomX}_{2} := (\LomX_{1} \times \LomX_{2}, \oqX, {\bf \jI}_{X}, \llX)$ and 
$\textbf{\LomY} = \textbf{\LomY}_{1} \times \textbf{\LomY}_{2} := (\LomY_{1} \times \LomY_{2}, \oqY, {\bf \jI}_{Y}, \llY)$.
\end{center} 
A map $\kap :\LomX_{1} \times \LomX_{2} \rightarrow \LomY_{1} \times \LomY_{2}$ is called a \textcolor{red}{\it direct product $O$-homomorphism} if it
is an $O$-homomorphism.
\end{definition}

\begin{theorem}\label{Tohdp-230524}
Let $\kap_{1} :\LomX_{1} \rightarrow \LomY_{1}$
 and $\kap_{2} :\LomX_{2} \rightarrow \LomY_{2}$ be $O$-homomorphisms from $\textbf{\LomX}_{1}$ to $\textbf{\LomY}_{1}$ and from  
		$\textbf{\LomX}_{2}$ to 	$\textbf{\LomY}_{2}$, 
respectively. Let $\textbf{\LomX} = \textbf{\LomX}_{1} \times \textbf{\LomX}_{2} := (\LomX_{1} \times \LomX_{2}, \oqX, {\bf \jI}_{X}, \llX)$ and 
$\textbf{\LomY} = \textbf{\LomY}_{1} \times \textbf{\LomY}_{2} := (\LomY_{1} \times \LomY_{2}, \oqY, {\bf \jI}_{Y}, \llY)$ be direct product OBCI-algebras. Then a map $\kap :\LomX_{1} \times \LomX_{2} \rightarrow \LomY_{1} \times \LomY_{2}$ such that 
	\begin{align}\label{hdp-230507}
		\kap : {\textbf{\LomX}} \rightarrow {\textbf{\LomY}}, ~(\ix, \kx) \mapsto (\kap_{1}(\ix), \kap_{2}(\kx))
	\end{align}
is a direct product $O$-homomorphism.
\end{theorem}

\begin{proof}
Suppose that a map $\kap :\LomX_{1} \times \LomX_{2} \rightarrow \LomY_{1} \times \LomY_{2}$ satisfies \eqref{hdp-230507}. We first verify the following homomorphism of $\kap$.
	\begin{align}\label{hdp}
		\kap((\ix, \kx) \oqX (\jxi, \jxk)) = \kap(\ix, \kx) \oqY \kap(\jxi, \jxk).
	\end{align}
We have
		$$\kap((\ix, \kx) \oqX (\jxi, \jxk)) = \kap(\ix \opA \jxi, \kx \opB \jxk)$$
by the definition of $\oqX$;
		$$\kap(\ix \opA \jxi, \kx \opB \jxk) = (\kap_{1}(\ix \opA \jxi), \kap_{2}(\kx \opB \jxk))$$
by \eqref{hdp-230507};
		$$(\kap_{1}(\ix \opA \jxi), \kap_{2}(\kx \opB \jxk)) = (\kap_{1}(\ix) \opC \kap_{1}(\jxi), \kap_{2}(\kx) \opD \kap_{2}(\jxk))$$
by \eqref{homo};
		$$(\kap_{1}(\ix) \opC \kap_{1}(\jxi), \kap_{2}(\kx) \opD \kap_{2}(\jxk)) = (\kap_{1}(\ix), \kap_{2}(\kx)) \oqY (\kap_{1}(\jxi), \kap_{2}(\jxk))$$
by the definition of $\oqY$; and
		$$(\kap_{1}(\ix), \kap_{2}(\kx)) \oqY (\kap_{1}(\jxi), \kap_{2}(\jxk)) = \kap(\ix, \kx) \oqY \kap(\jxi, \jxk)$$
by \eqref{hdp-230507}. Hence \eqref{hdp} holds true.
%%%%%%%%%%%%%%%%%%%%%%%%%%%%%%%%%%%%%%%%%%%%%%%%%%%%%%%%%%%%%%%%%%%%%%%%%%

We next verify the following $O$-mapping of $\kap$. 
	\begin{align}\label{ohdp}
		{\bf \jI}_{X} \llX (\ix, \kx) \oqX (\jxi, \jxk) ~\Rightarrow ~{\bf \jI}_{Y} \llY \kap(\ix, \kx) \oqY \kap(\jxi, \jxk).
	\end{align}
Suppose that ${\bf \jI}_{X} \llX (\ix, \kx) \oqX (\jxi, \jxk).$ Using \eqref{OBCI-5}, we have $(\ix, \kx) \llX (\jxi, \jxk)$, and so
\begin{center} 
$\jIA \leA \ix \opA \jxi$ and $\jIB \leB \kx \opB \jxk$
\end{center} 
 using \eqref{OBCI-5} and the definition of $\llX$. Then
\begin{center} 
$\jIC \leC \kap_{1}(\ix) \opC \kap_{1}(\jxi)$ and $\jID \leD \kap_{2}(\kx) \opD \kap_{2}(\jxk)$
\end{center} 
by \eqref{O-map}, and so
$$(\kap_{1}(\ix), \kap_{2}(\kx)) \llY (\kap_{1}(\jxi), \kap_{2}(\jxk))$$
 by \eqref{OBCI-5} and the definition of $\llY$. Hence
$${\bf \jI}_{Y} \llY (\kap_{1}(\ix), \kap_{2}(\kx)) \oqY (\kap_{1}(\jxi), \kap_{2}(\jxk))$$
 by \eqref{OBCI-5}, and so
$${\bf \jI}_{Y} \llY \kap(\ix, \kx) \oqY \kap(\jxi, \jxk)$$
by \eqref{hdp-230507}. This ensures that \eqref{ohdp} holds true.
%%%%%%%%%%%%%%%%%%%%%%%%%%%%%%%%%%%%%%%%%%%%%%%%%%%%%%%%%%%%%%%%%%%%%%%%%%

Therefore $\kap$ is a direct product $O$-homomorphism.
\end{proof}

\begin{definition}\label{def-kdp}
Suppose that $\kap_{1} :\LomX_{1} \rightarrow \LomY_{1}$
 and $\kap_{2} :\LomX_{2} \rightarrow \LomY_{2}$ are mappings from $\textbf{\LomX}_{1}$ to $\textbf{\LomY}_{1}$ and from $\textbf{\LomX}_{2}$ to $\textbf{\LomY}_{2}$, respectively; $\textbf{\LomX} = \textbf{\LomX}_{1} \times \textbf{\LomX}_{2} := (\LomX_{1} \times \LomX_{2}, \oqX, {\bf \jI}_{X}, \llX)$ and 
$\textbf{\LomY} = \textbf{\LomY}_{1} \times \textbf{\LomY}_{2} := (\LomY_{1} \times \LomY_{2}, \oqY, {\bf \jI}_{Y}, \llY)$ are direct product OBCI-algebras; and a map $\kap :\LomX_{1} \times \LomX_{2} \rightarrow \LomY_{1} \times \LomY_{2}$ satisfies \eqref{hdp-230507}. A subset $\LomA$ of $\LomX_{1} \times \LomX_{2}$ such that $\LomA = \LomA_{1} \times \LomA_{2}$, where $\LomA_{1} \subseteq \LomX_{1}$ and $\LomA_{2} \subseteq \LomX_{2}$, is called 
the \textcolor{red}{\it direct product kernel} of $\kap$, denoted by ${\rm ker}_{DP}(\kap)$, if both $\LomA_{1}$ and $\LomA_{2}$ satisfy $\eqref{ker}$.
\end{definition}

\begin{theorem}\label{Tkdp-230524}
Suppose that $\kap_{1} :\LomX_{1} \rightarrow \LomY_{1}$
 and $\kap_{2} :\LomX_{2} \rightarrow \LomY_{2}$ are mappings from $\textbf{\LomX}_{1}$ to $\textbf{\LomY}_{1}$ and from $\textbf{\LomX}_{2}$ to $\textbf{\LomY}_{2}$, respectively; $\textbf{\LomX} = \textbf{\LomX}_{1} \times \textbf{\LomX}_{2} := (\LomX_{1} \times \LomX_{2}, \oqX, {\bf \jI}_{X}, \llX)$ and 
$\textbf{\LomY} = \textbf{\LomY}_{1} \times \textbf{\LomY}_{2} := (\LomY_{1} \times \LomY_{2}, \oqY, {\bf \jI}_{Y}, \llY)$ are direct product OBCI-algebras; and a map $\kap :\LomX_{1} \times \LomX_{2} \rightarrow \LomY_{1} \times \LomY_{2}$ satisfies \eqref{hdp-230507}. Then ${\rm ker}(\kap_{1}) \times {\rm ker}(\kap_{2})$ is the direct product kernel of $\kap$.
\end{theorem}

\begin{proof}
Suppose that $\ix \in {\rm ker}(\kap_{1})$ and $\kx \in {\rm ker}(\kap_{2})$ for $\ix \in \LomX_{1}$ and $\kx \in \LomX_{2}$, and ${\bf \jIY} = (\jIC, \jID)$.  
We verify that
	\begin{align}\label{kdp}
		(\ix, \kx) \in {\rm ker}(\kap_{1}) \times {\rm ker}(\kap_{2}) ~\Leftrightarrow ~ {\bf \jIY} \llY \kap(\ix, \kx).
	\end{align}
We have
		$$(\ix, \kx) \in {\rm ker}(\kap_{1}) \times {\rm ker}(\kap_{2}) = (\ix \in {\rm ker}(\kap_{1}), \kx \in {\rm ker}(\kap_{2}))$$
by the Cartesian product of $\times$;
		$$(\ix \in {\rm ker}(\kap_{1}), \kx \in {\rm ker}(\kap_{2})) ~\Leftrightarrow ~ (\jIC \leC \kap_{1}(\ix), \jID \leD \kap_{2}(\kx))$$
by \eqref{ker};
		$$(\jIC \leC \kap_{1}(\ix), \jID \leD \kap_{2}(\kx)) ~\Leftrightarrow ~ (\jIC, \jID) \llY (\kap_{1}(\ix), \kap_{2}(\kx))$$
by the definition of $\llY$; and
		$$(\jIC, \jID) \llY (\kap_{1}(\ix), \kap_{2}(\kx)) ~\Leftrightarrow ~ {\bf \jIY} \llY \kap(\ix, \kx)$$
by \eqref{hdp-230507}. This ensures that \eqref{kdp} holds true. Therefore ${\rm ker}(\kap_{1}) \times {\rm ker}(\kap_{2})$ is the direct product kernel ${\rm ker}_{DP}(\kap)$ of $\kap$.
\end{proof}

\begin{corollary}\label{Ckdp-230524}
Suppose that $\kap_{1} :\LomX_{1} \rightarrow \LomY_{1}$
 and $\kap_{2} :\LomX_{2} \rightarrow \LomY_{2}$ are $O$-homomorphisms from $\textbf{\LomX}_{1}$ to $\textbf{\LomY}_{1}$ and from $\textbf{\LomX}_{2}$ to 	$\textbf{\LomY}_{2}$, respectively; $\textbf{\LomX} = \textbf{\LomX}_{1} \times \textbf{\LomX}_{2} := (\LomX_{1} \times \LomX_{2}, \oqX, {\bf \jI}_{X}, \llX)$ and 
$\textbf{\LomY} = \textbf{\LomY}_{1} \times \textbf{\LomY}_{2} := (\LomY_{1} \times \LomY_{2}, \oqY, {\bf \jI}_{Y}, \llY)$ are direct product OBCI-algebras; and a direct product $O$-homomorphism $\kap :\LomX_{1} \times \LomX_{2} \rightarrow \LomY_{1} \times \LomY_{2}$ satisfies \eqref{hdp-230507}. Then ${\rm ker}(\kap_{1}) \times {\rm ker}(\kap_{2})$ is the direct product kernel of $\kap$.
\end{corollary}

\begin{theorem}\label{T2kdp-230524}
Suppose that $\kap_{1} :\LomX_{1} \rightarrow \LomY_{1}$
 and $\kap_{2} :\LomX_{2} \rightarrow \LomY_{2}$ are mappings from $\textbf{\LomX}_{1}$ to $\textbf{\LomY}_{1}$ and from $\textbf{\LomX}_{2}$ to 	$\textbf{\LomY}_{2}$, respectively; $\textbf{\LomX} = \textbf{\LomX}_{1} \times \textbf{\LomX}_{2} := (\LomX_{1} \times \LomX_{2}, \oqX, {\bf \jI}_{X}, \llX)$ and 
$\textbf{\LomY} = \textbf{\LomY}_{1} \times \textbf{\LomY}_{2} := (\LomY_{1} \times \LomY_{2}, \oqY, {\bf \jI}_{Y}, \llY)$ are direct product OBCI-algebras; and a map $\kap :\LomX_{1} \times \LomX_{2} \rightarrow \LomY_{1} \times \LomY_{2}$ satisfies \eqref{hdp-230507}. We can express the direct product kernel ${\rm ker}_{DP}(\kap)$ of $\kap$ by  
${\rm ker}_{DP}(\kap) = {\rm ker}(\kap_{1}) \times {\rm ker}(\kap_{2})$ of $\textbf{\LomX}_{1}$ and $\textbf{\LomX}_{2}$, respectively.
\end{theorem}

\begin{proof}
Suppose that ${\rm ker}_{DP}(\kap)$ is the kernel of $\textbf{\LomX}_{1} \times \textbf{\LomX}_{2}$; $\pi_{\LomX_{1}}$ and $\pi_{\LomX_{2}}$ are projections such that $\pi_{\LomX_{1}} : \LomX_{1} \times \LomX_{2} \rightarrow \LomX_{1}$ and 
$\pi_{\LomX_{2}} : \LomX_{1} \times \LomX_{2} \rightarrow \LomX_{2}$, respectively; and ${\rm ker}(\kap_{1}):=\pi_{\LomX_{1}}({\rm ker}_{DP}(\kap))$ and ${\rm ker}(\kap_{2}):=\pi_{\LomX_{2}}({\rm ker}_{DP}(\kap))$. Let $\ix \in {\rm ker}(\kap_{1})$ and $\kx \in {\rm ker}(\kap_{2})$ for $\ix \in \LomX_{1}$ and $\kx \in \LomX_{2}$. We verify 
$$(\ix, \kx) \in {\rm ker}_{DP}(\kap) ~\Leftrightarrow ~\ix \in {\rm ker}(\kap_{1}), \kx \in {\rm ker}(\kap_{2}).$$
Note that ${\rm ker}_{DP}(\kap)$ is the kernel of $\LomX_{1} \times \LomX_{2}$. We obtain  
$$(\ix, \kx) \in {\rm ker}_{DP}(\kap) ~\Leftrightarrow ~ (\jIC, \jID) = {\bf \jI}_{Y} \llY \kap(\ix, \kx)$$
by \eqref{ker}, and so
$$(\jIC, \jID) \llY \kap(\ix, \kx) ~\Leftrightarrow ~(\jIC, \jID) \llY (\kap_{1}(\ix), \kap_{2}(\kx))$$
by \eqref{hdp-230507}. Then
$$(\jIC, \jID) \llY (\kap_{1}(\ix), \kap_{2}(\kx)) ~\Leftrightarrow ~\jIC \leC \kap_{1}(\ix), \jID \leD \kap_{2}(\kx)$$
by the definition of $\llY$, and so 
$$\jIC \leC \kap_{1}(\ix), \jID \leD \kap_{2}(\kx) ~\Leftrightarrow ~\ix \in {\rm ker}(\kap_{1}), \kx \in {\rm ker}(\kap_{2})$$
by \eqref{ker}.
Hence $(\ix, \kx) \in {\rm ker}_{DP}(\kap)$ if and only if $\ix \in {\rm ker}(\kap_{1})$ and $\kx \in {\rm ker}(\kap_{2})$. Therefore it holds true that ${\rm ker}_{DP}(\kap) = {\rm ker}(\kap_{1}) \times {\rm ker}(\kap_{2})$. 
\end{proof}

\begin{corollary}\label{C2kdp-230524}
Suppose that $\kap_{1} :\LomX_{1} \rightarrow \LomY_{1}$
 and $\kap_{2} :\LomX_{2} \rightarrow \LomY_{2}$ be $O$-homomorphisms from $\textbf{\LomX}_{1}$ to $\textbf{\LomY}_{1}$ and from $\textbf{\LomX}_{2}$ to 	$\textbf{\LomY}_{2}$, respectively; $\textbf{\LomX} = \textbf{\LomX}_{1} \times \textbf{\LomX}_{2} := (\LomX_{1} \times \LomX_{2}, \oqX, {\bf \jI}_{X}, \llX)$ and 
$\textbf{\LomY} = \textbf{\LomY}_{1} \times \textbf{\LomY}_{2} := (\LomY_{1} \times \LomY_{2}, \oqY, {\bf \jI}_{Y}, \llY)$ are direct product OBCI-algebras; and an $O$-homomorphism $\kap :\LomX_{1} \times \LomX_{2} \rightarrow \LomY_{1} \times \LomY_{2}$ satisfies \eqref{hdp-230507}.  We can express the direct product kernel ${\rm ker}_{DP}(\kap)$ of $\kap$ by 
${\rm ker}_{DP}(\kap) = {\rm ker}(\kap_{1}) \times {\rm ker}(\kap_{2})$ of $\textbf{\LomX}_{1}$ and $\textbf{\LomX}_{2}$, respectively.
\end{corollary}

The following are examples of Theorem \ref{T2kdp-230524} and Corollary \ref{C2kdp-230524}.

\begin{example}\label{ET2kdp1-230525}
Suppose that $\kap_{1} :\LomX_{1} \rightarrow \LomY_{1}$ are the $O$-map $\kap : \LomX \rightarrow \LomY$ in Example \ref{EqvO2h-230429} and $\kap_{2} :\LomX_{2} \rightarrow \LomY_{2}$ the $O$-homomorphism $\kap : \LomX \rightarrow \LomY$ in Example \ref{EqE161-230520}. Let $\textbf{\LomX} = \textbf{\LomX}_{1} \times \textbf{\LomX}_{2} := (\LomX_{1} \times \LomX_{2}, \oqX, {\bf \jI}_{X}, \llX)$ and 
$\textbf{\LomY} = \textbf{\LomY}_{1} \times \textbf{\LomY}_{2} := (\LomY_{1} \times \LomY_{2}, \oqY, {\bf \jI}_{Y}, \llY)$ be direct product OBCI-algebras, and a map $\kap :\LomX_{1} \times \LomX_{2} \rightarrow \LomY_{1} \times \LomY_{2}$ satisfy \eqref{hdp-230507}. Then \eqref{hdp} does not hold and so $\kap$ is not an $O$-homomorphism. 
To verify this, we show
	\begin{align*}
		\kap((\p, \jI) \oqX (\jI, \jI)) \neq \kap(\p, \jI) \oqY \kap(\jI, \jI).
	\end{align*}
The left part is as follows. 
	\begin{align*}
		& \indent \kap((\p, \jI) \oqX (\jI, \jI))
		\\&= \kap(\p \opA \jI, \jI \opB \jI)
		\\&= (\kap_{1}(\p \opA \jI), \kap_{2}(\jI \opB \jI))
		\\&= (\kap_{1}(\p), \kap_{2}(\jI))
		\\&= (\9, \jI).
	\end{align*}
The right part is as follows. 
	\begin{align*}
		& \indent \kap(\p, \jI) \oqY \kap(\jI, \jI)
		\\&= (\kap_{1}(\p), \kap_{2}(\jI)) \oqY (\kap_{1}(\jI), \kap_{2}(\jI))
		\\&= (\kap_{1}(\p) \opC \kap_{1}(\jI), \kap_{2}(\jI) \opD \kap_{2}(\jI))
		\\&= (\9 \opC \8, \jI \opD \jI)
		\\&= (\8, \jI).
	\end{align*}
Hence $(\9, \jI) \neq (\8, \jI),$ and so \eqref{hdp} does not hold.

${\rm ker}(\kap_{1})=\{\1, \jI\}$ and ${\rm ker}(\kap_{2})=\{\jI\}$. Then 
$${\rm ker}(\kap_{1}) \times {\rm ker}(\kap_{2})= \{(\ix, \kx)\in \LomX_{1} \times \LomX_{2} \mid \ix\in {\rm ker}(\kap_{1}), \, \kx\in {\rm ker}(\kap_{2})\}$$
is the kernel of $\textbf{\LomX}_{1} \times \textbf{\LomX}_{2} :=  (\LomX_{1} \times \LomX_{2}, \oqX, {\bf \jI}_{X}, \llX)$, where ${\bf \jI}_{X} = (\jI, \jI)$.
\end{example}

\begin{example}\label{ET2kdp2-230524}
Suppose that $\kap_{1} :\LomX_{1} \rightarrow \LomY_{1}$ are the homomorphism $\kap : \LomX \rightarrow \LomX$ in Example \ref{ETkcGV70-230520} and $\kap_{2} :\LomX_{2} \rightarrow \LomY_{2}$ the $O$-homomorphism $\kap : \LomX \rightarrow \LomY$ in Example \ref{EqE161-230520}. Let $\textbf{\LomX} = \textbf{\LomX}_{1} \times \textbf{\LomX}_{2} := (\LomX_{1} \times \LomX_{2}, \oqX, {\bf \jI}_{X}, \llX)$ and 
$\textbf{\LomY} = \textbf{\LomY}_{1} \times \textbf{\LomY}_{2} := (\LomY_{1} \times \LomY_{2}, \oqY, {\bf \jI}_{Y}, \llY)$ be direct product OBCI-algebras, and a map $\kap :\LomX_{1} \times \LomX_{2} \rightarrow \LomY_{1} \times \LomY_{2}$ satisfy \eqref{hdp-230507}. Then \eqref{ohdp} does not hold and so $\kap$ is not an $O$-homomorphism. 
To verify this, we show
	\begin{align*}
	&(\3, \jI) = {\bf \jI}_{X} \llX (\0, \jI) \oqX (\3, \jI),
	\\& (\3, \jI) = {\bf \jI}_{Y} \not\llY \kap(\0, \jI) \oqY \kap(\3, \jI).
	\end{align*}
The first case can be proved as follows: 
\begin{center}
$(\3, \jI) \llX (\0, \jI) \oqX (\3, \jI)$ if and only if\\
$\3 \leA \0 \opA \3 = \1$ and $\jI \leB \jI \opB \jI = \jI.$
\end{center}
The second case can be proved as follows:
\begin{center}
$\3 = \jIC \not\leC \kap_{1}(\0) \opC \kap_{1}(\3) = \1 \opC \3 = \0$ and\\
$\jI = \jID \leD \kap_{2}(\jI) \opD \kap_{2}(\jI) = \jI \opD \jI = \jI.$
\end{center}
Hence 
\begin{center}
${\bf \jI}_{X} \llX (\0, \jI) \oqX (\3, \jI) = (\1, \jI)$ but\\
${\bf \jI}_{Y} \not\llY \kap(\0, \jI) \oqY \kap(\3, \jI) = (\0, \jI),$
\end{center}
and so \eqref{ohdp} does not hold.

${\rm ker}(\kap_{1})=\{\1, \3\}$ and ${\rm ker}(\kap_{2})=\{\jI\}$. Then 
$${\rm ker}(\kap_{1}) \times {\rm ker}(\kap_{2})= \{(\ix, \kx)\in \LomX_{1} \times \LomX_{2} \mid \ix\in {\rm ker}(\kap_{1}), \, \kx\in {\rm ker}(\kap_{2})\}$$
is the kernel of $\textbf{\LomX}_{1} \times \textbf{\LomX}_{2} :=  (\LomX_{1} \times \LomX_{2}, \oqX, {\bf \jI}_{X}, \llX)$, where ${\bf \jI}_{X} = (\3, \jI)$.
\end{example}

\begin{example}\label{EC2kdp-230524}
Suppose that $\textbf{\LomX}_{1}$ is $\textbf{\LomX}$ in Example \ref{ETkcGV70-230520} and $\kap_{1} :\LomX_{1} \rightarrow \LomY_{1}$ is an automorphism as the identity map. Then it is clear that $\kap_{1}$ is an $O$-homomorphism. Suppose that $\kap_{2} :\LomX_{2} \rightarrow \LomY_{2}$ is the $O$-homomorphism $\kap : \LomX \rightarrow \LomY$ in Example \ref{EqE161-230520}. Let $\textbf{\LomX} = \textbf{\LomX}_{1} \times \textbf{\LomX}_{2} := (\LomX_{1} \times \LomX_{2}, \oqX, {\bf \jI}_{X}, \llX)$ and 
$\textbf{\LomY} = \textbf{\LomY}_{1} \times \textbf{\LomY}_{2} := (\LomY_{1} \times \LomY_{2}, \oqY, {\bf \jI}_{Y}, \llY)$ be direct product OBCI-algebras, and a map $\kap :\LomX_{1} \times \LomX_{2} \rightarrow \LomY_{1} \times \LomY_{2}$ satisfy \eqref{hdp-230507}. Then $\kap$ is an $O$-homomorphism by Theorem \ref{Tohdp-230524}. Note that ${\rm ker}(\kap_{1})=\{\1, \3\}$ and ${\rm ker}(\kap_{2})=\{\jI\}$. Then 
${\rm ker}(\kap_{1}) \times {\rm ker}(\kap_{2})$ is the kernel of $\textbf{\LomX}_{1} \times \textbf{\LomX}_{2}$
as in Example \ref{ET2kdp2-230524}.
\end{example}

\begin{theorem}\label{TkqP435-230524}
Suppose that $\kap_{1} :\LomX_{1} \rightarrow \LomY_{1}$
 and $\kap_{2} :\LomX_{2} \rightarrow \LomY_{2}$ are mappings from $\textbf{\LomX}_{1}$ to $\textbf{\LomY}_{1}$ and from $\textbf{\LomX}_{2}$ to 	$\textbf{\LomY}_{2}$, respectively; $\textbf{\LomX} = \textbf{\LomX}_{1} \times \textbf{\LomX}_{2} := (\LomX_{1} \times \LomX_{2}, \oqX, {\bf \jI}_{X}, \llX)$ and 
$\textbf{\LomY} = \textbf{\LomY}_{1} \times \textbf{\LomY}_{2} := (\LomY_{1} \times \LomY_{2}, \oqY, {\bf \jI}_{Y}, \llY)$ are direct product OBCI-algebras; and a map $\kap :\LomX_{1} \times \LomX_{2} \rightarrow \LomY_{1} \times \LomY_{2}$ satisfies \eqref{hdp-230507}. 
	For subsets $\LomK_{\LomX_{1}}$ and $\LomK_{\LomX_{2}}$ of $\textbf{\LomX}_{1}$ and 
	$\textbf{\LomX}_{2}$, respectively, we define two sets as follows:
$${\mathcal {\LomK}}^{\jID}:=\{(\ix, \kx)\in \LomX_{1} \times \LomX_{2}\mid 
     \ix\in \LomK_{\LomX_{1}}, \, \jID \leD \kap(\kx)\}.$$ 
$${\mathcal {\LomK}}^{\jIC}:=\{(\ix, \kx)\in \LomX_{1} \times \LomX_{2}\mid  
    \jIC \leC \kap(\ix), \, \kx\in \LomK_{\LomX_{2}}\}.$$
	If $\LomK_{\LomX_{1}}$ and $\LomK_{\LomX_{2}}$ are  kernels of 
	$\textbf{\LomX}_{1}$ and 
	$\textbf{\LomX}_{2}$, respectively, then 
${\mathcal {\LomK}}^{\jID}$ and ${\mathcal {\LomK}}^{\jIC}$  are  kernels of 
	$\textbf{\LomX}_{1} \times \textbf{\LomX}_{2}$.
\end{theorem}

\begin{proof}
Suppose that $\LomK_{\LomX_{1}}$ and $\LomK_{\LomX_{2}}$ are kernels of 
	$\textbf{\LomX}_{1}$ and 
	$\textbf{\LomX}_{2}$, respectively. Clearly ${\bf \jIX}=(\jIA, \jIB)\in {\mathcal \LomK}^{\jID}$.
Let  $(\ix, \kx) \in {\mathcal \LomK}^{\jID}$ for $(\ix, \kx) \in \LomX_{1} \times \LomX_{2}$. 
Then 
$\ix\in \LomK_{\LomX_{1}}$ and $\jID \leD \kap(\kx)$, and so $\kx\in \LomK_{\LomX_{2}}$ by \eqref{ker}. Hence ${\mathcal \LomK}^{\jID}$ is the kernel of 
$\textbf{\LomX}_{1} \times \textbf{\LomX}_{2}$.
In a similar manner, we can show that ${\mathcal \LomK}^{\jIC}$  is the kernel of 
$\textbf{\LomX}_{1} \times \textbf{\LomX}_{2}$.
\end{proof}

\begin{corollary}\label{CkqP435-230524}
Suppose that $\kap_{1} :\LomX_{1} \rightarrow \LomY_{1}$
 and $\kap_{2} :\LomX_{2} \rightarrow \LomY_{2}$ are $O$-homomorphisms from $\textbf{\LomX}_{1}$ to $\textbf{\LomY}_{1}$ and from $\textbf{\LomX}_{2}$ to 	$\textbf{\LomY}_{2}$, respectively; $\textbf{\LomX} = \textbf{\LomX}_{1} \times \textbf{\LomX}_{2} := (\LomX_{1} \times \LomX_{2}, \oqX, {\bf \jI}_{X}, \llX)$ and 
$\textbf{\LomY} = \textbf{\LomY}_{1} \times \textbf{\LomY}_{2} := (\LomY_{1} \times \LomY_{2}, \oqY, {\bf \jI}_{Y}, \llY)$ are direct product OBCI-algebras; and an $O$-homomorphism $\kap :\LomX_{1} \times \LomX_{2} \rightarrow \LomY_{1} \times \LomY_{2}$ satisfies \eqref{hdp-230507}. 
	For subsets $\LomK_{\LomX_{1}}$ and $\LomK_{\LomX_{2}}$ of $\textbf{\LomX}_{1}$ and 
	$\textbf{\LomX}_{2}$, respectively, we define two sets ${\mathcal {\LomK}}^{\jID}$ and ${\mathcal {\LomK}}^{\jIC}$ as in Theorem \ref{TkqP435-230524}.
	If $\LomK_{\LomX_{1}}$ and $\LomK_{\LomX_{2}}$ are  kernels of 
	$\textbf{\LomX}_{1}$ and 
	$\textbf{\LomX}_{2}$, respectively, then 
${\mathcal {\LomK}}^{\jID}$ and ${\mathcal {\LomK}}^{\jIC}$  are  kernels of 
	$\textbf{\LomX}_{1} \times \textbf{\LomX}_{2}$.
\end{corollary}

Finally we notice that ${\mathcal {\LomK}}^{\jID}$ and ${\mathcal {\LomK}}^{\jIC}$ in Theorem \ref{TkqP435-230524} and Corollary \ref{CkqP435-230524} are the same since $(\ix, \kx)\in {\mathcal {\LomK}}^{\jID}$ if and only if $(\ix, \kx)\in {\mathcal {\LomK}}^{\jIC}$, for $(\ix, \kx)\in \LomX_{1} \times \LomX_{2}$.

%%%%%%%%%%%%%%%%%%%%%%%%%%%%%%%%%%%%%%%%%%%%%%%%%%%%%%%%%%%%
\section{Conclusion}
%%%%%%%%%%%%%%%%%%%%%%%%%%%%%%%%%%%%%%%%%%%%%%%%%%%%%%%%%%%%

In this paper we introduced ordered homomorphisms of OBCI-algebras as order-preserving homomorphisms of OBCI-algebras and dealt with related properties. More exactly, first we defined the notions of the ordered homomorphism and kernel of OBCI-algebras. Next, we addressed properties related to (ordered) subalgebras, (ordered) filters and direct products of OBCI-algebras. 

We have some future works or open problems as follows. First, ordered homomorphisms and kernels of OBCI-algebras have to be studied in a more specific context. For example, ordered homomorphisms and kernels of OBCI-algebras related to (ordered) Y-filters, (ordered) R-filters and (ordered) J-filters need to be studied. Second, kernels of OBCI-algebras has to be generalized into ordered ones. As OBCI-algebras are a generalization of BCI-algebras with an underlying partial order, similarly we can introduce a generalization of the kernel of OBCI-algebras as ordered kernels of OBCI-algebras.

%%%%%%%%%%%%%%%%%%%%%%%%%%%%%%%%%%%%%%%%%%%%%%%%%%%%%%%%%%%%
\section{Declarations}
%%%%%%%%%%%%%%%%%%%%%%%%%%%%%%%%%%%%%%%%%%%%%%%%%%%%%%%%%%%%

%The authors declare that we do not have any competing financial interests or personal relationships that could have appeared to influence the work reported in this paper.

\subsection*{Data availability}
Data sharing not applicable to this article as datasets were neither generated nor analysed.

\subsection*{Compliance with ethical standards}
The authors declare that they have no conflict of interest.

 \end{document}